\documentclass[a4paper,12pt]{article}
\usepackage{pdfsync}
\usepackage{pgf,tikz}
\usepackage{amsfonts}
\usepackage{amssymb}
\usepackage{amsmath}
\usepackage{mathptmx}
\usepackage{cancel}
\usepackage{bbm}
\usepackage{pstricks}
\usepackage{psfrag}
\usepackage{mathrsfs}
\usepackage{graphicx}
\def\be{\begin{equation}}
\def\ee{\end{equation}}
\def\bea{\begin{eqnarray}}
\def\eea{\end{eqnarray}}
\def\beann{\begin{eqnarray*}}
\def\eeann{\end{eqnarray*}}
\definecolor{OliveGreen}{RGB}{40,220,40}

\def\ns{\hspace{-1mm}}
\newcommand{\real}{{\mathbb{R}}}

\def\spacingset#1{\def\baselinestretch{#1}\small\normalsize}
\spacingset{1.47}

\newtheorem{lemma}{Lemma}
\newtheorem{theorem}{Theorem}
\newtheorem{remark}{Remark}

\newtheorem{definition}{Definition}

\newtheorem{example}{Example}[section]
\newtheorem{assumption}{Assumption}[section]

\def\be{\begin{equation}}
\def\ee{\end{equation}}
\def\bea{\begin{eqnarray}}
\def\eea{\end{eqnarray}}
\def\beann{\begin{eqnarray*}}
\def\eeann{\end{eqnarray*}}

\def\ns{\hspace{-1mm}}

\def\proof{\noindent{\bf{\em Proof:}\ \ }}
\def\QED{\mbox{\rule[0pt]{1.5ex}{1.5ex}}}
\def\endproof{\hspace*{\fill}~\QED\par\endtrivlist\unskip}

\newcommand{\defi}{\stackrel{\text{\tiny def}}{=}}

\def\tra{\top}
\def\tp{^\top}

\definecolor{Royalblue}{cmyk}{1,0.30,0.2,0.2}

\newcommand{\complex}{{\mathbb{C}}}

\def\gR{{\cal R}}

\def\bmat{\left[ \begin{array}}
\def\emat{\end{array} \right]}

\def\bmat{\left[ \begin{array}}
\def\emat{\end{array} \right]}
\def\bsmat{\left[ \begin{smallmatrix}}
\def\esmat{\end{smallmatrix} \right]}

\def\gR{{\cal R}}

\def\i{{i}}

\begin{document}

\begin{titlepage}
\title{\vspace{-45mm}
Foundations of negative imaginary systems theory and relations with positive real systems\thanks{Partially supported by the Australian Research Council under the grant FT120100604 and by the
Italian Ministry for Education and Research (MIUR) under PRIN grant n. 20085FFJ2Z. Corresponding author L.~Ntogramatzidis. Tel. +61-8-92663143. }\vspace{10mm}}
\author{Augusto Ferrante$^\dagger$, Alexander Lanzon$^\ddagger$ and Lorenzo Ntogramatzidis$^\star$}
\date{$^\dagger${\small Dipartimento di Ingegneria dell'Informazione,\\[-2pt]
         Universit\`a di Padova, via Gradenigo, 6/B -- 35131 Padova, Italy. \\[-2pt]
         {\tt augusto@dei.unipd.it} 
   \\[8pt]
   $^\ddagger$ Control Systems Centre, School of Electrical and Electronic Engineering, University of Manchester, Sackville Street, Manchester M13 9PL, UK. \\[-2pt]
         {\tt Alexander.Lanzon@manchester.ac.uk.} 
   \\[8pt]
           $^\star$Department of Mathematics and Statistics,\\[-2pt]
         Curtin University, Perth (WA), Australia. \\[-2pt]
         {\tt L.Ntogramatzidis@curtin.edu.au}       }
        }%
\thispagestyle{empty} \maketitle \thispagestyle{empty}
\begin{abstract}%
In this paper we lay the foundations of a not necessarily rational negative imaginary systems theory and its relations with positive real systems theory and, hence, with passivity.
In analogy with the theory of positive real functions, in our general framework, negative imaginary systems are defined in terms of a domain of analyticity of the transfer function and of a sign condition that must be satisfied in such domain.
In this way, on the one hand, our theory does not require to restrict the attention to systems with rational transfer function and, on the other hand | just by suitably selecting the domain of analyticity to be either the right half complex plane or the complement of the unit disc in the complex plane | we particularize our theory to both continuous-time and to discrete-time systems. Indeed, to the best of our knowledge, this is first time that discrete-time negative imaginary systems are studied in the literature.
In this work, we also aim to provide a unitary view of the different notions that have appeared so far in the literature within the framework of positive real and in the more recent theory of negative imaginary systems, and to show how these notions are characterized and linked to each other.
 A stability analysis result for the interconnection of discrete-time systems is also derived.\end{abstract}
%
%
\thispagestyle{empty}
\end{titlepage}

\section{Introduction}
\label{secintro}
The theory of positive real systems is one of the fundamental cornerstones of systems and control theory, and in particular of passivity theory. Given the extensive amount of contributions in this area, dating back from the early 1930s \cite{Brune-31}, it would be impossible to quote all of the relevant references. We consequently
refer the readers to two important monographs \cite{Andersen-V-73}, \cite{Brogliato-LME-07} for a summary of the historic and recent contributions on this problem.
A promising new development in the theory of dissipative systems theory has been the introduction of the notion of {\em negative imaginary systems}, see \cite{Lanzon-P-08,Xiong-PL-10,Mabrok-KPL-14-1} and the references cited therein.
With respect to positive realness, the definition of negative imaginary systems imposes a weaker restriction on the relative degree of the transfer function and does not prohibit all unstable zeros. Negative imaginary systems theory was found to be very suitable in a range of applications including modelling and control of undamped or lightly damped flexible structures with colocated position sensors and force actuators \cite{Petersen-Lanzon-10,BMP-12}, in nano-positioning control due to piezoelectric transducers and capacitive sensors (e.g. \cite{BM-09,MMB-11,Mabrok-KPL-14-2}) and in multi-agent networked systems (e.g. \cite{CH-10,Wang-LP-14}). The notion of negative imaginary systems specializes also to the important subclass of lossless negative imaginary systems \cite{Xiong-PL-12}.

In spite of the wealth of results that in just a few years have been presented and published on negative imaginary systems including extensions to infinite dimensional systems \cite{Opmeer-11}, Hamiltonian systems \cite{Schaft-11}, descriptor systems \cite{Mabrok-KPL-12} and mixtures of negative imaginary and small-gain properties \cite{Patra-L-11} to mention only a few, an important gap in the current literature -- that the present paper attempts to fill -- is the lack of a definition of negative imaginary (and strictly negative imaginary) function for discrete-time systems. Furthermore, so far \cite{Ferrante-N-13} has been to the best of the authors' knowledge the only contribution which attempted to address the general case of a definition of negative imaginary system for non necessarily rational transfer functions, and then recovered the standard definition given in the foundational paper \cite{Lanzon-P-08} for the symmetric rational case.
 However, several aspects of the core theory of negative imaginary systems remained unexplored in \cite{Ferrante-N-13}. For example, the notion of {\em strictly negative imaginary system} has never been defined in the general case of a non-rational transfer function. This remaining gap will also be filled in this paper as it is essential in studying stability interconnections of both rational and non rational negative imaginary systems.



 To summarize, the main contribution of this paper is to present a general and foundational perspective of the recent theory of negative imaginary systems, and their relation with the classical theory of positive real systems. As a byproduct, we fill some important gaps that have so far remained open. In particular, 
\begin{enumerate}

\item As pointed out in \cite{Brogliato-LME-07}, since the early studies in the 1960s, there has been a proliferation of definitions of various types of strictly positive real systems. Our first aim is to follow the approach of \cite{Brogliato-LME-07} in the attempt of defining
different notions of strictly negative imaginary system and establishing a parallel between these definitions and their positive real counterparts. The standard notion of strictly negative imaginary system introduced in the literature so far corresponds to only one of these definitions. We will define, examine and characterize other notions of strictly negative imaginary functions.
\item The notion of discrete-time negative imaginary systems is introduced for the first time. This definition is given in the general non-rational setting and then is specialized for rational transfer functions, and expressed in terms of a sign constraint on the unit circle. We also introduce different notions of {\em strictly negative imaginary} discrete-time transfer functions.
 Finally, the relations between discrete-time and continuous-time negative imaginary systems are elucidated. We also provide a discrete-time negative imaginary lemma which yields a complete state-space characterization of discrete-time negative imaginary systems and a stability analysis result for the interconnection of discrete-time negative imaginary systems.
\end{enumerate}

\noindent{\bf Notation.} Given a matrix $A$, the symbol $A^\tra$ denotes the transpose of $A$ and $A^\ast$ denotes the complex conjugate transpose of $A$. 
{ We denote by $\sigma(A)$ the set of singular values of the matrix $A$ and with $\min \sigma(A)$ the smallest of such singular values.}
Recall that given a real rational function $G(s)$ and a simple pole $p\in \complex$ of $G(s)$, we have a unique decomposition $G(s)=G_1(s) +A/(s-p)$, where $G_1(s)$ is a rational function which is analytic in an open set containing $p$ and the (non-zero) matrix $A$ is the residue corresponding to the pole $p$. If $p$ is a double pole of $G(s)$, we have the unique decomposition $G(s)=G_1(s) +A_1/(s-p)+ A_2/(s-p)^2$, where the matrix $A_1$ is the residue corresponding to the pole $p$. In this case, by analogy, we define the (non-zero) matrix $A_2$ to be the {\em quadratic residue} corresponding to the pole $p$.
If $G(s)$ has a pole at infinity, it can be uniquely decomposed as
$G(s)=G_1(s) +P(s)$, where $G_1(s)$ is a rational proper function and $P(s)=\sum_{i=1}^k A_is^i$
is a homogeneous polynomial in $s$.
We refer to $A_i$ as the $i$-th coefficient in the expansion at infinity of $G(s)$.
The usual notations of $\ge 0$ and $>0$ are used to denote positive
semidefiniteness and positive definiteness of Hermitian matrices, respectively. Let $G : \complex \longrightarrow {\complex^{m \times m}}$ be analytic or harmonic in a certain region $\Omega$ of $\complex$; then $G$ is said to have full normal rank if there exists $s\in\Omega$ such that $\det[G(s)]\neq 0$.
Given complex matrices $S_1, S_2$ and complex vectors $y_1, y_2, u_1, u_2, \alpha, \beta$ of compatible dimension satisfying
$ \bsmat y_1 \\[1mm] \alpha \esmat=S_1 \bsmat u_1\\[1mm] \beta\esmat$ and $\bsmat \beta \\[1mm] y_2\esmat=S_2 \bsmat \alpha \\[1mm] u_2\esmat$,
let $S_1\star S_2$ denote the Redheffer star product which maps $\bsmat u_1 \\[1mm] u_2\esmat$ to $\bsmat y_1\\[1mm] y_2\esmat$.
Furthermore, $F_l(S_1,S_2^{(1,1)})$ (resp.~$F_u(S_2,S_1^{(2,2)})$) denote the lower (resp.~upper) linear fractional transformation.
Let $[P,Q]$ denote the positive feedback interconnection between systems $P$ and $Q$. Let $\bar{\lambda}(A)$ denote the largest eigenvalue of a square matrix $A$ that has only real eigenvalues.

\section{The Continuous-Time Case}
In this section, for the sake of completeness we briefly recall the most important notions and results of positive real and negative imaginary systems for the continuous-time case.

\begin{definition}
Let $F : \complex \longrightarrow {\complex^{m \times m}}$ be a continuous-time transfer function. Then, $F(s)$ is {\em continuous positive real} (C-PR) if
\begin{itemize}
\item $F(s)$ is analytic in $\{s\in\complex:\mathfrak{Re}\{s\} > 0\}$;
\item $F(s)$ is real when $s$ is real and positive;
\item $F(s)+F(s)^\ast \ge 0$ for all $s\in\complex$ such that $\mathfrak{Re}\{s\}> 0$.
\end{itemize}
\end{definition}

\begin{lemma}\label{CNSPR}
Let $F : \complex \longrightarrow {\complex^{m \times m}}$ be a continuous-time, real, rational transfer function. Then, $F(s)$ is C-PR if and only if
\begin{itemize}
\item $F(s)$ has no poles in $\{s\in\complex:\mathfrak{Re}\{s\} > 0\}$;
\item $F(i\,\omega)+F(i\,\omega)^\ast\ge 0$ for all $\omega\in\real$ such that $s=i\,\omega$ is not a pole of $F(s)$;
\item if $i\,\omega_{\scriptscriptstyle 0}$ is a pole of any element of $F(s)$, it is a simple pole with Hermitian and positive semidefinite residue. In particular, if $\omega_{\scriptscriptstyle 0}$ is finite, the residue is
 \[
K_{\scriptscriptstyle 0} \defi \lim_{s\to i \omega_{\scriptscriptstyle 0}} (s-i\,\omega_{\scriptscriptstyle 0})\,F(s),
\]
while if $\omega_{\scriptscriptstyle 0}$ is infinite, the residue is 
 \[
K_\infty \defi \lim_{\omega \to \infty} \frac{F(i\,\omega)}{i\,\omega}.
\]
\end{itemize}
\end{lemma}

{ We now present our definitions of strictly positive real systems.
We warn the reader that many different definitions have been proposed for this concept that can indeed be distinguished via several grades of strength, see e.g. \cite{Brogliato-LME-07,Khalil-02}.
In this paper, we shall only need two of such grades | that will be referred to as
{\em strongly} and {\em weakly} strictly positive realness | and we only briefly hint to a third, extra-strong,
grade.

\begin{definition}\label{C-SSPR}
Let $F : \complex \longrightarrow {\complex^{m \times m}}$ be a continuous-time, real transfer function. Then, $F(s)$ is {\em continuous strongly strictly positive real} (C-SSPR) if
for some $\epsilon>0$, the transfer function $F(s-\epsilon)$ is C-PR and $F(s)+F(-s)^\top$ has full normal rank.
\end{definition}

}

Now, we show that C-SSPR as defined in Definition~\ref{C-SSPR} can be equivalently checked via a strict sign condition in the domain of analyticity.

\begin{lemma}
\label{C-SSPRonAllRHP}
Let $F : \complex \longrightarrow {\complex^{m \times m}}$ be a continuous-time, real transfer function. Then, $F(s)$ is C-SSPR if and only if there exists $\epsilon>0$ such that
\begin{description}
\item{\bf {\em (i)}\;} $F(s)$ is analytic in $\{s\in\complex : \mathfrak{Re}\{s\}>-\epsilon\}$;
\item{\bf {\em (ii)}\;} $F(s)+F(s)^\ast > 0$ for all $s \in \{s \in \complex:\, \mathfrak{Re}\{s\}>-\epsilon \}$.
\end{description}
\end{lemma}
The proof of this result can be carried out by adapting the proof of Lemma \ref{C-SSNIonAllRHP} in the sequel, and it is therefore omitted.

{ The following result, see \cite[Theorem 2.47]{Brogliato-LME-07} and \cite[Lemma 6.1]{Khalil-02},
shows that in the case of rational functions the property of C-SSPR is equivalent to an analyticity condition and a sign condition restricted to the extended imaginary axis.}
\begin{theorem}\label{cns-c-sspr}
\label{brogliato}
Let $F : \complex \longrightarrow {\complex^{m \times m}}$ be a continuous-time, real, rational, proper transfer function. Then $F(s)$ is C-SSPR if and only if
\begin{enumerate}
\item $F(s)$ has all its poles with strictly negative real parts;
\item $F(i\omega)+F(-i\omega)^\top>0$ for all $\omega \in \real$;
\item one of the three conditions is satisfied:\footnote{ We write this property as three separate conditions to elucidate all the possible situations that may occur. It is clear, however, that the third condition is the more general and encompasses the first and the second. Essentially what this condition says is that the smallest singular value of $F(i\omega)+F(i\omega)^\top$ cannot tend to zero faster than $1/\omega^2$.}
\begin{itemize}
\item $F(\infty)+F(\infty)^\top>0$
\item $F(\infty)+F(\infty)^\top=0$ and $\lim_{\omega \to \infty} \omega^2 [F(i\omega)+F^\top(-i\omega)]>0$
\item $F(\infty)+F(\infty)^\top\ge0$ but not zero nor non-singular, and there exist $\sigma>0$ and $\delta>0$ such that
\[
\min \sigma \left[\omega^2 \,\left(F(i\omega)+F(-i\omega)^\top\right)\right]\ge \sigma_{\scriptscriptstyle 0}, \quad \forall |\omega|\ge \delta.
\]
\end{itemize}
\end{enumerate}
\end{theorem}

{ In some situations the concept of C-SSPR is too restrictive: indeed in the rational case where there are finitely many poles and zeros, it is useful to introduce the following weaker definition.
\begin{definition}\label{C-WSPR}
Let $F : \complex \longrightarrow {\complex^{m \times m}}$ be a continuous-time, real, rational, proper transfer function. Then, $F(s)$ is {\em continuous weakly strictly positive real} (C-WSPR) if
the first two properties of Theorem \ref{cns-c-sspr} hold.
\end{definition}
}

\begin{remark}
{\em
As shown in \cite[pp.~238--240]{Khalil-02}, examples of C-SSPR transfer functions are $F(s)=\frac{1}{s+a}$ (with $a>0$).
{ Notice that, if in the definition of C-SSPR we removed the full normal rank condition on $F(s)+F(-s)^\top$, we would have that functions such as
 $F(s)=\frac{1}{s+1}\bsmat 1 && 1 \\[1mm] 1 && 1 \esmat$ would be C-SSPR which is unacceptable as, in this case, a result like Theorem \ref{brogliato} would not hold.
}
An example of a transfer function which is C-PR but not C-SSPR is the following:
\[
F(s)=\frac{s+3}{(s+1)(s+2)}.
\]
Indeed, 1) in Theorem \ref{brogliato} is satisfied. Moreover, given $\varepsilon>0$, a simple calculation gives
\bea
\label{formulaesempio}
F(i\,\omega-\varepsilon)+F(i\,\omega-\varepsilon)^\ast = 2\,\frac{6+6\,\varepsilon^2-\varepsilon^3-\varepsilon\,(11+\omega^2)}{\left[ \omega^2+(2-\varepsilon)^2\right]
\left[ \omega^2+(1-\varepsilon)^2\right]},
\eea
which is strictly positive on the imaginary axis (i.e., when $\varepsilon=0$), so that 2) in Theorem \ref{brogliato} also holds. On the other hand, 3) in Theorem \ref{brogliato} is not satisfied. Indeed, in this case $F(\infty)+F(\infty)^\top=0$, but $\lim_{\omega \to \infty} \omega^2 [F(i\omega)+F^\top(-i\omega)]= \lim_{\omega \to \infty} \frac{12\,\omega^2}{(\omega^2+4)(\omega^2+1)}=0$. This result is consistent with Definition \ref{C-SSPR}. In fact, (\ref{formulaesempio}) shows that for any arbitrarily small $\varepsilon>0$, by taking a sufficiently large $\omega>0$, the numerator of $F(i\,\omega-\varepsilon)+F(i\,\omega-\varepsilon)^\ast$ can be rendered negative. In other words, $F(i\,\omega)+F(-i\,\omega)^\top$ is positive definite for all $\omega>0$, but no matter how small we choose $\varepsilon>0$, if $\omega>0$ is sufficiently large we can find $F(i\,\omega-\varepsilon)+F(i\,\omega-\varepsilon)^\ast<0$, and therefore $F(s-\varepsilon)$ is not C-PR for any $\varepsilon>0$.
Finally, we recall that in \cite{Brogliato-LME-07} also an ``extra strong" form of strict positive realness is defined which essentially correspond to coercivity of the corresponding spectral density.
}

\end{remark}

We now introduce the following standing assumption, that will be used throughout the rest of the paper.

\begin{assumption}
We henceforth restrict our attention to only symmetric transfer functions.
\end{assumption}

{As discussed in \cite{Ferrante-N-13}, the case of symmetric transfer function is the most important and interesting one, because it encompasses both the scalar case, and the case of a transfer function of a reciprocal $m$-port electrical network.\footnote{The only way to obtain a non-symmetric transfer function of an $m$-port electrical network is to employ gyrators, whose physical implementation requires the use of active components but that cannot be physically implemented with arbitrary precision.} Moreover, to the best of the authors' knowledge, all the negative imaginary transfer functions considered or studied in the literature so far are symmetric (see e.g. the transfer functions from a force actuator to a corresponding collocated position sensor | for instance, a piezoelectric sensor | in a lightly damped or undamped structure), even though the real, rational definitions of negative imaginary systems in \cite{Lanzon-P-08,Xiong-PL-10,Mabrok-KPL-14-1} allow for non-symmetric transfer functions.}

\begin{definition}
\label{defus}
Let $G : \complex \longrightarrow {\complex^{m \times m}}$ be a {\em continuous-time} transfer function. Then, $G$ is {\em continuous negative imaginary} (C-NI) if
\begin{description}
\item{ \hspace{-0.2cm}{\bf {\em (i)}\;}} $G(s)$ is analytic in $\{s\in\complex : \mathfrak{Re}\{s\}>0\}$;
\item{ \hspace{-0.2cm}{\bf {\em (i)}\;}} $i\,[ G(s)-G(s)^\ast ] \ge 0$ for all $s \in \{s \in \complex:\, \mathfrak{Re}\{s\}>0 \;\;\text{and}\;\;\mathfrak{Im}\{s\} >0 \}$;
\item{ \hspace{-0.2cm}{\bf {\em (iii)}\;}} $i\,[ G(s)-G(s)^\ast ] = 0$ for all $s \in \{s \in \complex:\, \mathfrak{Re}\{s\}>0 \;\;\text{and}\;\;\mathfrak{Im}\{s\} =0 \}$;
\item{ \hspace{-0.2cm}{\bf {\em (iv)}\;}} $i\,[ G(s)-G(s)^\ast ] \le 0$ for all $s \in \{s \in \complex:\, \mathfrak{Re}\{s\}>0 \;\;\text{and}\;\;\mathfrak{Im}\{s\} <0 \}$.
\end{description}
\end{definition}
The following result, which was proven in \cite{Ferrante-N-13}, provides a characterisation of rational NI systems in terms of a domain of analyticity and conditions referred to the imaginary axis.
\begin{lemma}
\label{lemmarestriction}
Let $G : \complex \longrightarrow {\complex^{m \times m}}$ be a continuous-time, real, rational transfer function. Then $G(s)$ is C-NI if and only if
\begin{description}
\item{ \hspace{-0.2cm}{\bf {\em (i)}\;}} $G(s)$ has no poles in $\mathfrak{Re}\{s\}>0$;
\item{ \hspace{-0.2cm}{\bf {\em (ii)}\;}} $i\,[G(i\, \omega)-G(i \omega)^\ast]\ge 0$ for all $\omega \in (0,\infty)$ except for the values of $\omega$ where $i\omega$ is a pole of $G(s)$;
\item{ \hspace{-0.2cm}{\bf {\em (iii)}\;}} {if $s = i\,\omega_{\scriptscriptstyle 0}$, with $\omega_{\scriptscriptstyle 0} \in (0,\infty)$, is a pole of $G(s)$, then it is a simple
pole and the corresponding {\em residual matrix}\footnote{Notice that $K_{\scriptscriptstyle 0}$ is the product of the imaginary unit $i$ by the residue in $\omega_{\scriptscriptstyle 0}$.} $K_{\scriptscriptstyle 0}= \lim_{s \to i\,\omega_{\scriptscriptstyle 0}} (s-i\,\omega_{\scriptscriptstyle 0})\,i\,G(s)$ is Hermitian and
 positive semidefinite;
\item{ \hspace{-0.2cm}{\bf {\em (iv)}\;}} if $s = 0$ is a pole of $G(s)$, then it is at most a double pole.
 Moreover, both its residual and its quadratic residual (when present) are positive semidefinite Hermitian matrices;
\item{ \hspace{-0.2cm}{\bf {\em (v)}\;}} if $s=\infty$ is a pole of $G(s)$, then it is at most a double pole.
 Moreover, both the coefficients in the expansion at infinity of $G(s)$ are negative semidefinite Hermitian matrices.}
 \end{description}
 \end{lemma}

\begin{remark}
{\em
We observe that $\frac{1}{s}$ and $\frac{1}{s^2}$ are negative imaginary, whereas $-\frac{1}{s^2}$ is not. Note that when there are poles on the imaginary axis, the $D$-contour is indented infinitesimally to the right and hence the Nyquist plot changes phase rapidly at large magnitudes around the frequency of the pole(s) on the imaginary axis. From the complete Nyquist plot it is evident that $\frac{1}{s}$ and $\frac{1}{s^2}$ are negative imaginary, but $-\frac{1}{s^2}$ is not.
}
\end{remark}

We recall the following important result, which established a relationship between C-PR and C-NI transfer functions, see \cite{Lanzon-P-08,Xiong-PL-10,Ferrante-N-13}.
\begin{theorem}
\label{prni}
 Let $G(s)$ be a real, rational, proper, symmetric negative imaginary transfer function matrix. Then $F(s)\defi s[G(s)-G(\infty)]$ is positive real.
Conversely, let $F(s)$ be real, rational, symmetric positive real transfer function matrix. Then $G(s)\defi (1/s)F(s) +D$ is symmetric negative imaginary for any symmetric matrix $D$.
\end{theorem}

We now adapt the definition of {strongly strictly} positive real function to the negative imaginary case.

{
\begin{definition}
\label{C-SSNI}
Let $G : \complex \longrightarrow {\complex^{m \times m}}$ be a continuous-time, real transfer function. Then, $G(s)$ is {\em continuous strongly strictly negative imaginary} (C-SSNI) if
for some $\epsilon>0$, the transfer function $G(s-\epsilon)$ is C-NI and $i\,[ G(s)-G(-s)^\top ]$ has full normal rank.
\end{definition}

\begin{remark}
{\em Note that the full normal rank condition is essential in the above definition as this class of systems will be needed for internal stability of positive feedback interconnections of C-NI and C-SSNI systems. If we were not to impose the full normal rank condition on the C-SSNI class, then the feedback interconnection of a C-NI system and a C-SSNI system would not be internally stable as demonstrated via the following simple example: Let $P(s)=\bsmat 1 && 1 \\[1mm] 1 && 1\esmat$ which is clearly C-NI and let $Q(s)=\frac{1}{s+1}\bsmat 1 && 1 \\[1mm] 1 && 1\esmat$ which fulfils all properties of C-SSNI except for the full normal rank condition. The positive feedback interconnection of $P(s)$ and $Q(s)$ is not internally stable as there exists a closed-loop pole at $s=3$.
}
\end{remark}

Now, we show that C-SSNI as defined in Definition~\ref{C-SSNI} can be equivalently checked via
conditions on the imaginary axis. To this aim, we need some preliminary results, starting with writing necessary and sufficient conditions on the domain of analyticy for a system to be C-SSNI.
\begin{lemma}
\label{C-SSNIonAllRHP}
Let $G : \complex \longrightarrow {\complex^{m \times m}}$ be a continuous-time, real transfer function. Then, $G(s)$ is C-SSNI if and only if there exists $\epsilon>0$ such that
\begin{description}
\item{ \hspace{-0.2cm}{\bf {\em (i)}\;}} $G(s)$ is analytic in $\{s\in\complex : \mathfrak{Re}\{s\}>-\epsilon\}$;
\item{ \hspace{-0.2cm}{\bf {\em (ii)}\;}} $i\,[ G(s)-G(s)^\ast ] > 0$ for all $s \in \{s \in \complex:\, \mathfrak{Re}\{s\}>-\epsilon \;\;\text{and}\;\;\mathfrak{Im}\{s\} >0 \}$;
\item{ \hspace{-0.2cm}{\bf {\em (iii)}\;}} $i\,[ G(s)-G(s)^\ast ] = 0$ for all $s \in \{s \in \complex:\, \mathfrak{Re}\{s\}>-\epsilon \;\;\text{and}\;\;\mathfrak{Im}\{s\} =0 \}$;
\item{ \hspace{-0.2cm}{\bf {\em (iv)}\;}} $i\,[ G(s)-G(s)^\ast ] < 0$ for all $s \in \{s \in \complex:\, \mathfrak{Re}\{s\}>-\epsilon \;\;\text{and}\;\;\mathfrak{Im}\{s\} <0 \}$.
\end{description}
\end{lemma}

\proof
Definition~\ref{C-SSNI} trivially gives equivalence to the existence of $\epsilon>0$ such that conditions {\bf {\em (i)}-{\em (iv)}} are satisfied with non-strict inequalities in {\bf {\em (ii)}} and {\bf {\em (iv)}} on $i\,[ G(s)-G(s)^\ast ]$. We hence only need to show that if $G$ is C-SSNI, then the inequalities in {\bf {\em (ii)}} and {\bf {\em (iv)}} are indeed strict. We prove only that {\bf {\em (ii)}} is strict since {\bf {\em (iv)}} follows by symmetry.
Let $G$ be analytic in $\complex_{-\varepsilon} \defi\{s\in\complex : \mathfrak{Re}\{s\}>-\epsilon\}$ and assume by contradiction that there exist $s_{\scriptscriptstyle 0}\in \{s \in \complex:\, \mathfrak{Re}\{s\}>-\epsilon \;\;\text{and}\;\;\mathfrak{Im}\{s\} >0 \}$ and a nonzero vector $v$ such that $v^\ast (i\,[ G(s_{\scriptscriptstyle 0})-G(s_{\scriptscriptstyle 0})^\ast ])v =0$.
Let $\varepsilon_1<\epsilon$ be such that $\mathfrak{Re}\{s_{\scriptscriptstyle 0}\}>-\varepsilon_1$.
Since $G$ is analytic in $\complex_{-\varepsilon}$, $v^\ast (i\,[ G(s)-G(s)^\ast ])v$
is harmonic in the same domain so that, by considering an arbitrarily large real number $M$ and the compact set ${\mathcal C} \defi\{s\in\complex : M\geq\mathfrak{Re}\{s\}\geq -\varepsilon_1 \;\;\text{and}\;\; M\geq\mathfrak{Im}\{s\}\geq 0 \}\subset \complex_{-\varepsilon}$, if $v^\ast (i\,[ G(s)-G(s)^\ast ])v$ restricted to ${\mathcal C}$ attains its minimum at a point $s_{\scriptscriptstyle 0}$ in the interior of ${\mathcal C}$, then $v^\ast (i\,[ G(s)-G(s)^\ast ])v$ is constant. Clearly, $v^\ast (i\,[ G(s)-G(s)^\ast ])v\geq 0$ for all $s\in{\mathcal C}$ and, by taking $M$ sufficiently large, $s_{\scriptscriptstyle 0}$ is in the interior of ${\mathcal C}$ so that $v^\ast (i\,[ G(s)-G(s)^\ast ])v$
is constantly equal to $0$. This is a contradiction, since Definition~\ref{C-SSNI} requires that $i\,[ G(s)-G(-s)^\top ]$ has full normal rank.
\endproof

{
\begin{lemma}
\label{C-SSNI-pre-pre-pre-theorem}
Let $g : \complex \longrightarrow {\complex}$ be a scalar, continuous-time, real, rational, strictly proper transfer function. Assume that $g(s)$ is a C-NI function. Then, the relative degree of $g(s)$ is at most $2$ and all the finite zeros of $g(s)$ are in the closed left half-plane.
Moreover, if $i\,[g(i\omega)-g(i\omega)^\ast]>0$ for all $\omega \in (0,\infty)$, then all the finite zeros of $g(s)$ are in the open left half-plane.
\end{lemma}
\proof
As a consequence of \cite[Theorem 3.1]{Ferrante-N-13} we have that $f(s) \defi s\,g(s)$ is C-PR. Then, the relative degree of $f(s)$ is at most $1$ and all the finite zeros of $f(s)$ are in the closed left half-plane.
Therefore, the relative degree of $g(s)$ is at most $2$ and all the finite zeros of $g(s)$ are in the closed left half-plane.
Moreover, if $i\,[g(i\omega)-g(i\omega)^\ast]>0$, the only point of the imaginary axis in which $g$ could vanish is $0$. If, however, $g(0)=0$ then $f(s)$ would have a double zero at the origin which is in contrast with positive realness.
\endproof

{
\begin{remark}
{\em
Note that the strictly proper assumption in Lemma~\ref{C-SSNI-pre-pre-pre-theorem} is essential to this observation. Indeed, it is possible to have bi-proper transfer functions such as $g(s)=\frac{1-s}{1+s}$ that have zeros in the open right half-plane and are still C-SSNI and hence also C-NI. This is a crucial difference between PR functions (that are necessarily minimum phase) and NI functions.
}
\end{remark}
}

\begin{lemma}\label{C-SSNI-pre-pre-theorem}
Let $g : \complex \longrightarrow {\complex}$ be a scalar, continuous-time, real, rational, proper transfer function. Assume that $g(s)$ is a C-SSNI function.
 If $g(0)=0$, then the multiplicity of the zero in the origin of $g(s)$ is equal to $1$.\end{lemma}
\proof
Since $g(s)$ is a C-SSNI function, it has no poles in zero and we can expand $g(s)$ at the origin as
$$
g(s)=\sum_{k=h}^\infty r_k s^k,
$$
where $h$ is the multiplicity of the zero at the origin of $g$.
Let $s=\varepsilon e^{i\theta}$, $0<\theta<\pi$.
If $\varepsilon$ is sufficiently small,
$
i[g(s)-g(s)^\ast] $
has the same sign of
$ -2 r_h \varepsilon^h\sin(h\theta)$,
so that it can be positive for any $\theta\in(0,\pi)$ only if $h=1$.
%
\endproof
}

{
{ We now present necessary and sufficient conditions on the imaginary axis for a system to be C-SSNI.}
\begin{theorem}\label{C-SSNItheorem}
Let $G : \complex \longrightarrow {\complex^{m \times m}}$ be a continuous-time, real, rational, proper transfer function. Then $G(s)$ is C-SSNI if and only if
\begin{description}
\item{ \hspace{-0.2cm}{\bf {\em (i)}\;}} $G(s)$ has all its poles with strict negative real parts;
\item{ \hspace{-0.2cm}{\bf {\em (ii)}\;}} $i\,[G(i\omega)-G(i\omega)^\ast]>0$ for all $\omega \in (0,\infty)$;
\item{ \hspace{-0.2cm}{\bf {\em (iii)}\;}} There exist $\sigma_{\scriptscriptstyle 0}>0$ and $\delta>0$ such that
\begin{equation}\label{max-rel-deg}
\min\sigma[\omega^3 i\,[G(i\omega)-G(i\omega)^\ast]]>\sigma_{\scriptscriptstyle 0}\ \forall\ \omega \ge \delta ;
\end{equation}
\item{ \hspace{-0.2cm}{\bf {\em (iv)}\;}}
\begin{equation}\label{max-zero-origin}
Q \defi\lim_{\omega\rightarrow 0^+} (1/\omega) i\,[G(i\omega)-G(i\omega)^\ast]>0.
\end{equation}
\end{description}
\end{theorem}
}}

\proof
Necessity of {\bf {\em (i)}} and {\bf {\em (ii)}} is trivial from Lemma~\ref{C-SSNIonAllRHP}.
We now show necessity of condition {\bf {\em (iii)}}.
Essentially, we need to show that for any vector $v$ the relative degree of $i[g'(i\omega)-g'(i\omega)^\ast]$, where $g'(s) \defi v\tp G(s) v$, is at most $3$.
Assume by contradiction that this is not the case so that
$g(s) \defi g'(s)-g'(\infty)$ is a rational strictly proper C-SSNI function such that
$i[g(i\omega)-g(i\omega)^\ast]$ tends to zero, as $\omega\rightarrow\infty$, faster than $1/\omega^3$.
Then, it is easy to check that the relative degree of $g$ is at least $2$ and, in view of Lemma \ref{C-SSNI-pre-pre-pre-theorem}, the relative degree of $g$ is exactly $2$.
In view of Lemma \ref{C-SSNI-pre-pre-pre-theorem} we can write $g(s)$ as
$$
g(s)=K\frac{n(s)}{d(s)}=K\frac{s^{n-2}+a_{n-3}s^{n-3}+\dots +a_{\scriptscriptstyle 0}}{s^{n}+b_{n-1}s^{n-1}+\dots +b_{\scriptscriptstyle 0}},
$$
with $a_i$ and $b_i$ strictly positive.
By imposing that $i[g(i\omega)-g(i\omega)^\ast]$ tends to zero, as $\omega\rightarrow\infty$, faster than $1/\omega^3$, we get that $n\geq 3$ and $a_{n-3}=b_{n-1}$.
Now, we can compute
$$i[g(i\omega-\varepsilon)-g(i\omega-\varepsilon)^\ast]
=\frac{-4K\varepsilon\omega}{|d(i\omega-\varepsilon)|^2}
[(\varepsilon^2+\omega^2)^{n-2}+T_{2n-6}]
$$
with $T_{2n-6}$ being a polynomial in $\omega$ of degree equal to $2n-6$.
Therefore for a sufficiently large $\omega$, $i[g(i\omega)-g(i\omega)^\ast]$ is negative for any positive $\varepsilon$.

{
We now show necessity of condition {\bf {\em (iv)}}.
Assume that $G$ is C-SSNI. Then clearly the limit $Q$ defined in (\ref{max-zero-origin}) exists and is positive semi-definite.
Assume by contradiction that $Q$ is singular and let $v\in\ker Q$.
Let $g'(s) \defi v\tp G(s) v$. Clearly, $g'(s)$ is a rational proper C-SSNI function and
$g(s) \defi g'(s)-g'(\infty)$ is a rational strictly proper C-SSNI function such that
\begin{equation}\label{max-zero-origin-1}
\lim_{\omega\rightarrow 0} (1/\omega) i\,[g(i\omega)-g(i\omega)^\ast]=0.
\end{equation}
In view of Lemma \ref{C-SSNI-pre-pre-pre-theorem} we can write $g(s)$ as
$$
g(s)=K\frac{n(s)}{d(s)}=K\frac{1+a_1s+a_2s^2+\dots +a_ms^m}{1+b_1s+b_2s^2+\dots +b_ns^n},\quad m<n
$$
with $a_i$ and $b_i$ strictly positive.
Then (\ref{max-zero-origin-1}) implies $a_1=b_1$.
Notice now that
$$
g(s)-K=K\frac{n(s)-d(s)}{d(s)}$$
is C-SSNI as well so that the multiplicity of its zero in the origin is at most equal to $1$.
Therefore $a_1\neq b_1$.

{
As for sufficiency,
assume that $G(s)$ is real symmetric and rational and that it satisfies {\bf {\em (i)}}, {\bf {\em (ii)}}, {\bf {\em (iii)}} and {\bf {\em (iv)}}.
We now show that we can choose $\varepsilon>0$ in such a way that
\begin{equation}\label{tesi-intermedia}
i\,[G(-\varepsilon+i\omega)-G(-\varepsilon+i\omega)^\ast]>0,\ \forall\ \omega \in (0,\infty).
\end{equation}
In view of condition {\bf {\em (ii)}}, we have that for all $\omega_2>\omega_1>0$, there exists
$\varepsilon>0$ such that
\begin{equation}\label{tesi-intermedia-1}
i\,[G(-\varepsilon+i\omega)-G(-\varepsilon+i\omega)^\ast]>0,\ \forall\ \omega \in [\omega_1,\omega_2],
\end{equation}
so that
it is sufficient to show that given an arbitrarily small $\omega_1$ and an arbitrarily large $\omega_2$, there exists
$\varepsilon>0$ such that
\begin{equation}\label{tesi-intermedia-2}
i\,[G(-\varepsilon+i\omega)-G(-\varepsilon+i\omega)^\ast]>0,\ \forall\ \omega \in (0,\omega_1)
\end{equation}
and
\begin{equation}\label{tesi-intermedia-3}
i\,[G(-\varepsilon+i\omega)-G(-\varepsilon+i\omega)^\ast]>0,\ \forall\ \omega \in (\omega_2,\infty).
\end{equation}

As for (\ref{tesi-intermedia-2}), let $\delta \defi i\omega-\varepsilon$ and consider the following expansion of $G(\delta)$:
$$
G(\delta)=D_{\scriptscriptstyle 0}+\delta D_1+\delta^2 D_2+\dots
$$
which clearly converges for $\delta$ sufficiently small (if we considered a minimal realization $G(s)=C(sI-A)^{-1}B+D$, we would have $D_{\scriptscriptstyle 0} \defi D-CA^{-1}B$ and $D_i \defi-CA^{-i-1}B$, for $i>1$).
Notice that since $G(s)$ is real symmetric by standing assumption, $D_i=D_i\tp$.
Moreover, $Q \defi\lim_{\omega\rightarrow 0^+} (1/\omega) i\,[G(i\omega)-G(i\omega)^\ast]=-2D_1$, so that by assumption {\bf {\em (iv)}}, we have $D_1<0$.
Now a direct calculation gives
$$
i\,[G(-\varepsilon+i\omega)-G(-\varepsilon+i\omega)^\ast]=-\omega \,2\, D_1 + i\,\sum_{j=2}^\infty [\delta^j-(\delta^\ast)^j]\,D_j.$$
Now we observe that
$$i\,\sum_{j=2}^\infty [\delta^j-(\delta^\ast)^j]\,D_j=
-2\,\omega \,\sum_{j=2}^\infty \sum_{k=0}^{j-1} [\delta^k(\delta^\ast)^{j-1-k}]\,D_j,$$
so that
$$\| i\,\sum_{j=3}^\infty [\delta^j-(\delta^\ast)^j]\,D_j\|\leq
2\,\omega \,\sum_{j=2}^\infty j \varepsilon^{j-1}\|\,D_j\|=
2\,\omega \,\varepsilon \sum_{j=2}^\infty j \varepsilon^{j-2}\|D_j\|\leq
2\,\omega \,\sigma \varepsilon$$
for a certain $\sigma$ which does not increase as $\varepsilon$ tends to zero.
Since, by choosing a sufficiently small $\varepsilon$ we can make
$-D_1>\sigma\varepsilon I$, we have (\ref{tesi-intermedia-2}).

Now we prove (\ref{tesi-intermedia-3}).
Let $G(s)=C(sI-A)^{-1}B+D$ be a a minimal realization so that
$$
G(-\varepsilon+i\omega)=G(i\omega)+\varepsilon\Delta(i\omega)
$$
with $\Delta(s) \defi C(sI-A)^{-1}(sI-\varepsilon I-A)^{-1}B+D$.
We can expand $\Delta$ around infinity
as
$$
\Delta(i\omega)=\frac{CB}{(i\omega)^2}+\frac{\Delta_3}{(i\omega)^3} + \frac{\Delta_r(i\omega)}{(i\omega)^4}
$$
where $\Delta_3$ remains bounded as $\varepsilon$ tends to zero and $\Delta_r(i\omega)$ remains bounded as $\varepsilon$ tends to zero and $\omega$ tends to $+\infty$.
Then, we have
\beann 
i\,[G(-\varepsilon+i\omega)-G(-\varepsilon+i\omega)^\ast] \ns&\ns = \ns&\ns 
i\,[G(i\omega)-G(i\omega)^\ast]
+\frac{\varepsilon(-\Delta_3-\Delta_3\tp)}{\omega^3}\\
 \ns&\ns \ns&\ns
+\frac{i}{\omega^4}[\Delta_r(i\omega)-\Delta_r(i\omega)^\ast]
\eeann
so that, in view of condition {\bf {\em (iii)}}, (\ref{tesi-intermedia-3}) holds.

Now we can apply Lemma \ref{lemmarestriction} to the function $G(s-\varepsilon)$ and we immediately see that it is C-NI so that $G$ is C-SSNI
 }
\endproof

References \cite{SLPP-SCL-12,LSPP-CIS-11} and earlier define strictly negative imaginary systems by imposing only conditions {\bf {\em (i)}} and {\bf {\em (ii)}} of Theorem~\ref{C-SSNItheorem}.

We then introduce the following definition.
\begin{definition}
\label{defweak}
The continuous-time, real, rational, proper transfer function $G: \complex \longrightarrow \complex^{m \times m}$ is {\em continuous weakly strictly negative imaginary} (C-WSNI) if it satisfies conditions {\bf {\em (i)}} and {\bf {\em (ii)}} of Theorem~\ref{C-SSNItheorem}.
\end{definition}
}

The following two examples show that conditions {\bf {\em (iii)}} and {\bf {\em (iv)}} in Theorem \ref{C-SSNItheorem} are not implied by the first two, i.e., the notion of C-WSNI is indeed a weaker notion than that of C-SSNI.

\begin{example}
{\em
Consider the transfer function
\[
G(s)=\frac{2\,s+1}{(s+1)^2}.
\]
It is easily seen that $G(s)$ is C-NI. A simple calculation shows that
\bea
\label{fores1}
i\,\left[ G(i\,\omega-\varepsilon)-G(i\,\omega-\varepsilon)^\ast \right]=\frac{4\,\omega}{\left[ \omega^2+(1-\varepsilon)^2\right]^2}(\omega^2-\varepsilon+\varepsilon^2),
\eea
which proves that conditions {\bf {\em (i)}}, {\bf {\em (ii)}} and {\bf {\em (iii)}} in Theorem \ref{C-SSNItheorem} are satisfied; in particular, this means that $G(s)$ is C-WSNI. However, it is not C-SSNI, because in this case (\ref{max-zero-origin}) yields $Q=
\lim_{\omega\rightarrow 0^+} \frac{4\,\omega^2}{(\omega^2+1)^2}=0$. This result is consistent with Definition \ref{C-SSNI}. Indeed, for any $\varepsilon>0$, there always exists a sufficiently small $\omega>0$ such that the numerator in (\ref{fores1}) is negative.
}
\end{example}

\begin{example}
{\em
Consider the transfer function
\[
G(s)=\frac{s+3}{(s+1)^3}.
\]
Again, $G(s)$ is C-NI, and in this case
\bea
\label{fores2}
i\,\left[ G(i\,\omega-\varepsilon)-G(i\,\omega-\varepsilon)^\ast \right]=\frac{4\,\omega\,\left[ 4+6\,\varepsilon^2-\varepsilon^3-\varepsilon\,(\omega^2+9)\right]}{(1+\varepsilon^2-2\,\varepsilon+\omega^2)^3}.
\eea
Thus, conditions {\bf {\em (i)}}, {\bf {\em (ii)}} in Theorem \ref{C-SSNItheorem} are satisfied, which means that $G(s)$ is C-WSNI. Condition {\bf {\em (iv)}} in Theorem \ref{C-SSNItheorem} is also satisfied, since in this case (\ref{max-zero-origin}) gives $Q=
\lim_{\omega\rightarrow 0^+} \frac{16}{(\omega^2+1)^3}=16>0$. However, $G(s)$ is not C-SSNI because {\bf {\em (iii)}} in Theorem \ref{C-SSNItheorem} is not satisfied. Again, this result is consistent with Definition \ref{C-SSNI}, since for any $\varepsilon>0$, there always exists a sufficiently large $\omega>0$ such that the numerator in (\ref{fores2}) becomes negative.
}
\end{example}

The next lemma shows that the definition of C-WSNI corresponds to a sign property on the closed right-half plane.
\begin{lemma}
\label{C-WSNIonAllRHP}
Let $G : \complex \longrightarrow {\complex^{m \times m}}$ be a continuous-time, real, rational, proper transfer function. Then, $G(s)$ is C-WSNI if and only if there exists $\epsilon>0$ such that
\begin{description}
\item{ \hspace{-0.2cm}{\bf {\em (i)}\;}} $G(s)$ is analytic in $\{s\in\complex : \mathfrak{Re}\{s\}>-\epsilon\}$;
\item{ \hspace{-0.2cm}{\bf {\em (ii)}\;}} $i\,[ G(s)-G(s)^\ast ] > 0$ for all $s \in \{s \in \complex:\, \mathfrak{Re}\{s\}\geq 0 \;\;\text{and}\;\;\mathfrak{Im}\{s\} >0 \}$;
\item{ \hspace{-0.2cm}{\bf {\em (iii)}\;}} $i\,[ G(s)-G(s)^\ast ] = 0$ for all $s \in \{s \in \complex:\, \mathfrak{Re}\{s\}\geq 0 \;\;\text{and}\;\;\mathfrak{Im}\{s\} =0 \}$;
\item{ \hspace{-0.2cm}{\bf {\em (iv)}\;}} $i\,[ G(s)-G(s)^\ast ] < 0$ for all $s \in \{s \in \complex:\, \mathfrak{Re}\{s\}\geq 0 \;\;\text{and}\;\;\mathfrak{Im}\{s\} <0 \}$.
\end{description}
\end{lemma}
\proof
Sufficiency is trivial by restricting on the imaginary axis. Necessity can be proven as follows: if $G$ is C-WSNI, then {\bf {\em (i)}} is satisfied and $G$ is C-NI (from Lemma \ref{lemmarestriction}). Moreover, if $G$ is C-NI then {\bf {\em (ii)}}-{\bf {\em (iv)}} in Definition \ref{defus} are satisfied. Appending the imaginary axis properties of $G$ to the conditions {\bf {\em (ii)}}-{\bf {\em (iv)}} in Definition \ref{defus} (since $G$ is C-WSNI) yields {\bf {\em (ii)}}-{\bf {\em (iv)}} since $G$ fulfills {\bf {\em (i)}}.
\endproof

The following result shows that a relationship can be established between C-SSPR and C-SSNI transfer functions that is the counterpart of
Theorem \ref{prni}.

\begin{lemma}
Let $F: \complex \longrightarrow {\complex^{m \times m}}$ be a symmetric, real, rational, C-SSPR transfer function. Then, there exists $\varepsilon >0$ such that $G(s) \defi \frac{F(s)}{s+\varepsilon}+D$ is C-SSNI for any symmetric matrix $D$. Conversely, let $G: \complex \longrightarrow {\complex^{m \times m}}$ be a symmetric, real, rational C-SSNI transfer function. Then, there exists $\varepsilon >0$ such that $(s+\varepsilon)\left(G(s)-G(\infty)\right)$ is C-SSPR.
\end{lemma}
\proof
We start proving the first statement. Since $F(s)$ is C-SSPR, a value $\varepsilon >0$ exists such that $F(s-\varepsilon)$ is C-PR. Then, by Theorem \ref{prni} it is found that $\hat{G}(s) \defi \frac{F(s-\varepsilon)}{s}$ is C-NI. On the other hand, this implies that $G(s)=\frac{F(s)}{s+\varepsilon}$ is C-SSNI.\\
We prove the second statement. Since $G(s)$ is C-SSNI, there exists $\varepsilon >0$ such that $G(s-\varepsilon)$ is C-NI. Thus, by Theorem \ref{prni} we find that $s\,\left[ G(s-\varepsilon)-G(\infty)\right]$ is C-PR. This in turn implies that $(s+\varepsilon)\,\left[ G(s)-G(\infty)\right]$ is C-SSPR.
\endproof

\section{The Discrete-Time Case}
%

\subsection{Discrete-Time Positive Real Systems}
The definition of {\em discrete-time positive real} function was introduced for the first time in the literature by Hitz and Anderson in \cite{Hitz-A-69}, and is recalled below.

\begin{definition}{\sc [Hitz and Anderson, \cite{Hitz-A-69}]} \\
\label{defdspr}
The function $F: \complex \longrightarrow \complex^{m \times m}$ is {\em discrete positive real} (D-PR) if
\begin{itemize}
\item $F(z)$ is analytic in $\{z\in\complex:|z| >1\}$;
\item $F(z)$ is real when $z$ is real and positive;
\item $F(z)^\ast+F(z)\ge 0$ for all $|z| >1$.
\end{itemize}
\end{definition}

Similarly to what happens in the continuous-time,
for rational functions, discrete positive realness can be characterised in terms of conditions involving properties of the restriction of the matrix function to the unit circle.

\begin{theorem}{\sc \cite{Premaratne-J-94,Tao-I-90}}.
\label{restriction}
Let $F : \complex \longrightarrow {\complex^{m \times m}}$ be a discrete-time, real, rational, proper transfer function. 
Then, $F(z)$ is D-PR if and only if
\begin{itemize}
\item $F(z)$ is analytic in $\{z\in\complex:|z| > 1\}$;
\item $F(e^{i\,\theta})^\ast+F(e^{i \,\theta})\ge 0$ for all $\theta \in [0,2\pi)$ except for the values of $\theta$ for which $z=e^{i\,\theta}$ is a pole of $F(z)$;
\item If $z_{\scriptscriptstyle 0}=e^{i\,\theta_{\scriptscriptstyle 0}}$, with $\theta_{\scriptscriptstyle 0} \in [0,2\pi)$, is a pole of $F(z)$, then it is a simple pole and the normalized residual matrix
\beann
K_{\scriptscriptstyle 0} \defi \frac{1}{z_0}\,\lim_{z\to z_{\scriptscriptstyle 0}} (z-z_{\scriptscriptstyle 0}) \,F(z)
\eeann
 is Hermitian and positive semidefinite.
\end{itemize}
\end{theorem}

We now introduce the notion of discrete-time strongly strictly positive realness for discrete-time transfer functions.

\begin{definition}\label{D-SSPR}
Let $F : \complex \longrightarrow {\complex^{m \times m}}$ be a discrete-time, real, proper transfer function. Then, $F(z)$ is {\em discrete strongly strictly positive real} (D-SSPR) if,
for some $\delta\in (0,1)$, the transfer function $F(\delta z)$ is D-PR and $F(z)+F(1/z)^\top$ has full normal rank.
\end{definition}

The following simple result is the discrete-time counterpart of Theorem \ref{brogliato}.
\begin{theorem}\label{proponcondfordssprini}
Let $F : \complex \longrightarrow {\complex^{m \times m}}$ be a discrete-time, real, rational, proper transfer function.
Then, $F(z)$ is D-SSPR if and only if
\begin{itemize}
\item $F(z)$ has all its poles in a disc of radius $\rho \in [0,1)$;
\item $F(e^{i\,\theta})+F(e^{i\,\theta})^\ast>0$ for all $\theta \in [0,2\,\pi)$.
\end{itemize}
\end{theorem}
\proof
Necessity of the first condition is obvious. Necessity of the second immediately follows from the fact the unit circle is in the interior of the domain of analyticity and by the full normal rank assumption.
As for sufficiency,
since the unit circle is closed, condition $F(e^{i\,\theta})+F(e^{i\,\theta})^\ast>0$ for all $\theta \in [0,2\,\pi)$ implies coercivity, i.e. there exists $\sigma_{\scriptscriptstyle 0}>0$ such that
$F(e^{i\,\theta})+F(e^{i\,\theta})^\ast>\sigma_{\scriptscriptstyle 0}\,I$ for all $\theta \in [0,2\,\pi)$. Therefore, there exists $\rho \in [0,1)$ such that
$F(\rho e^{i\,\theta})+F(\rho e^{i\,\theta})^\ast>0$ for all $\theta \in [0,2\,\pi)$, so that
$F_1(z)\defi F(\rho z)$ is D-PR.
\endproof

{
\begin{remark}
{\em
The conditions of Theorem \ref{proponcondfordssprini} are much simpler than those of Theorem \ref{brogliato} because in view of the closure of the unit circle $\mathbb{T}$ (as opposed to the fact that the imaginary axis is not closed) positivity in $\mathbb{T}$ implies coercivity.
As we shall see later this is not the case for the NI systems for which the relevant boundary curve is the intersection between $\mathbb{T}$ and the open upper half complex plane.
Therefore, the relevant boundary curve is not closed as the zero and infinity discrete frequencies are not in this curve.
}
\end{remark}
}

The next result is the discrete-time counterpart of the so-called {\em positive real lemma}, a cornerstone of modern control theory that has generated an endless stream of literature.

\begin{lemma}{\sc [Discrete-Time Positive Real Lemma, \cite[Lemma 3]{Hitz-A-69}]} \\
\label{theoremAnderson}
Let $F : \complex \longrightarrow {\complex^{m \times m}}$ be a discrete-time, real, rational, proper transfer function with no poles in $|z|>1$ and simple poles only on $|z|=1$. Let $(A,B,C,D)$ be a minimal  realization of $F(z)$. Then $F(z)$ is discrete positive real if and only if there exists a real matrix $X=X^\top>0$ and real matrices $L$ and $W$ such that
\bea
&& X-A^\top\,X\,A=L^\top\,L, \label{prl1}\\
&& C^\top -A^\top\,X\,B=L^\top\,W, \label{prl2} \\
&& D^\top+D-B^\top\,X\,B=W^\top\,W. \label{prl3}
\eea
\end{lemma}

\subsection{Discrete-Time Negative Imaginary Functions}
We now present a definition of negative imaginary functions in the discrete-time case.
\begin{definition}
\label{dni}
Let $G: \complex \longrightarrow \complex^{m \times m}$ be a discrete-time, real transfer function. We say that $G(z)$ is \emph{discrete
negative imaginary} (D-NI) if
\begin{description}
\item{ \hspace{-0.2cm}{\bf {\em (i)}\;}} $G(z)$ is analytic in $\{z\in \complex\,:\,\,|z|>1\}$;
\item{ \hspace{-0.2cm}{\bf {\em (ii)}\;}} $i\,\left[G(z)-G(z)^\ast\right]\ge 0$ for all $z\in\{z\in \complex\,:\,\,|z|>1 \; \text{and} \;\, \mathfrak{Im}(z)>0\}$;
\item{ \hspace{-0.2cm}{\bf {\em (iii)}\;}} $i\,\left[G(z)-G(z)^\ast\right]= 0$ for all $z\in\{z\in \complex\,:\,\,|z|>1 \; \text{and} \;\, \mathfrak{Im}(z)=0\}$;
\item{ \hspace{-0.2cm}{\bf {\em (iv)}\;}} $i\,\left[G(z)-G(z)^\ast\right]\le 0$ for all $z\in\{z\in \complex\,:\,\,|z|>1 \; \text{and} \;\, \mathfrak{Im}(z)<0\}$.
\end{description}
\end{definition}

The conditions {\bf {\em (ii)}}-{\bf {\em (iv)}} in Definition \ref{dni} are a {\em skew imaginary condition} on the open set $\Omega=\{z\in \complex \,:\;|z|>1\}$.

\begin{remark}
{\em
If the real transfer function $G: \complex \longrightarrow \complex^{m \times m}$ satisfies the conditions in Definition \ref{dni}, then $G(z)$ is symmetric, i.e., $G(z)=G(z)^\top$ for all $z\in \complex$ such that $|z|>1$. This can be seen as follows: since $G(z)$ is real, if $z \in \real$ then $G(z)\in \real$. Let $z \in \real$ and $|z|>1$. From {\bf {\em (iii)}}, we get $G(z)=G(z)^\top$. Since this holds for all $z \in \real$ and $|z|>1$, the identity theorem of analytic functions ensures that this holds for all $z \in \complex$ in the domain of analyticity, i.e., $|z|>1$.
}
\end{remark}

Conditions {\bf {\em (iii)}}-{\bf {\em (iv)}} in Definition \ref{dni} are redundant in the real rational case, as the following result establishes.

\begin{lemma}
Let $G: \complex \longrightarrow \complex^{m \times m}$ be a discrete-time, real, rational transfer function. If $G(z)$ satisfies {\bf {\em (i)}}-{\bf {\em (ii)}} of Definition \ref{dni}, then it also satisfies {\bf {\em (iii)}}-{\bf {\em (iv)}}.
\end{lemma}
\proof If $G(z)$ satisfies {\bf {\em (ii)}}, then 
$i\,\left[G^\top(z)-G(\overline{z}) \right] \ge 0$ for all $z\in \complex$ such that $|z|>1$ and $\mathfrak{Im}(z)>0$, since $G({z})^\ast=G(\overline{z})^\top$.
 Defining $w \defi \overline{z}$, such condition can be re-written as
$i\,\left[G(w)^\ast-G(w) \right] \ge 0$ for all $w \in \complex$
 such that $|w|>1$ and $\mathfrak{Im}(w)<0$, which is exactly {\bf {\em (iv)}} of Definition \ref{dni}.
Finally, since {\bf {\em (ii)}} and {\bf {\em (iv)}} hold, then {\bf {\em (iii)}} must also hold by continuity.
\endproof

We now prove the counterpart of Theorem \ref{restriction} for the case of discrete-time symmetric negative imaginary functions.

\begin{lemma}
\label{resd}
Let $G: \complex \longrightarrow \complex^{m \times m}$ be a discrete-time, real, rational, proper transfer function. Then, $G(z)$ is D-NI if and only if
\begin{description}
\item{ \hspace{-0.2cm}{\bf {\em (i)}\;}} $G(z)$ has no poles in $|z|>1$;
\item{ \hspace{-0.2cm}{\bf {\em (ii)}\;}} $i\,[G(e^{i\theta})-G(e^{i\theta})^\ast]\ge 0$ for all $\theta \in (0,\pi)$ except for the values of $\theta$ for which $z=e^{i\,\theta}$ is a pole of $G(z)$;
\item{ \hspace{-0.2cm}{\bf {\em (iii)}\;}} {if $z_{\scriptscriptstyle 0} = e^{i\,\theta_{\scriptscriptstyle 0}}$, with $\theta_{\scriptscriptstyle 0} \in (0,\pi)$, is a pole of $G(z)$, then it is a simple
pole and the normalized residual matrix\footnote{Notice that $K_{\scriptscriptstyle 0}$ is the product of the imaginary unit $i$ by the residue in $z_{\scriptscriptstyle 0}$.} 
\bea
\label{normalresidu}
K_{\scriptscriptstyle 0}= \frac{1}{z_0}\,\lim_{z \to z_{\scriptscriptstyle 0}} (z-z_{\scriptscriptstyle 0})\,i\,G(z)
\eea
 is Hermitian and
 positive semidefinite;
\item{ \hspace{-0.2cm}{\bf {\em (iv)}\;}} if $z_{\scriptscriptstyle 0} = 1$ is a pole of $G(z)$, then it is at most a double pole.
 Moreover, its residual $A_1$ and its quadratic residual $A_2$ (when the pole is simple it is assumed that $A_2=0$)
are Hermitian matrices satisfying $A_2\geq 0$ and $A_1\geq A_2$;
\item{ \hspace{-0.2cm}{\bf {\em (v)}\;}} if $z_{\scriptscriptstyle 0} = -1$ is a pole of $G(z)$, then it is at most a double pole.
 Moreover, its residual $A_1$ and its quadratic residual $A_2$ (when the pole is simple it is assumed that $A_2=0$)
are Hermitian matrices satisfying $A_2\leq 0$ and $A_1\geq -A_2$.}
 \end{description}
 \end{lemma}


\proof
Since $G(z)$ is discrete-time real, symmetric and rational, define
$$G_c(s) \defi G\left(\frac{1+s}{1-s}\right).$$
Consider the identity
\[
z=\frac{1+s}{1-s},
\]
and let $z=\sigma+i\,\omega$. It is found that
\bea
\label{use}
s=\frac{z-1}{z+1}=\frac{\sigma^2+\omega^2-1}{(\sigma+1)^2+\omega^2}+2\,\i\,\frac{\omega}{(\sigma+1)^2+\omega^2}.
\eea
The following facts are easy to check:
\begin{enumerate}
\item
$G(z)$ is D-NI if and only if
$G_c(s)$ is C-NI.
Indeed, in view of (\ref{use}), $G(z)$ is analytic in $|z|>1$ if and only if $G_c(s)$ is analytic in $\mathfrak{Re}\{s\}>0$.
The rest of the proof of this part follows directly from the definitions, using the fact that $\mathfrak{Im}\{z\}>0$ (resp. $\mathfrak{Im}\{z\}<0$ and $\mathfrak{Im}\{z\}=0$) is equivalent to $\omega>0$ (resp. $\omega<0$ and $\omega=0$), which in turn is equivalent to $\mathfrak{Im}\{s\}>0$ (resp. $\mathfrak{Im}\{s\}<0$ and $\mathfrak{Im}\{s\}=0$).
\item
\begin{description}
\item{ \hspace{-0cm} {\bf {\em (i)}}} $G(z)$ has no poles in $|z|>1$ if and only if $G_c(s)$ has no poles in $\mathfrak{Re}\{s\}>0$;
\item{ \hspace{-0cm} {\bf {\em (ii)}}}
Let $z_{\scriptscriptstyle 0} \defi e^{i\,\theta_{\scriptscriptstyle 0}}$ with $\theta_{\scriptscriptstyle 0} \in (0,\pi)$. Using (\ref{use}) we see that
\[
s_{\scriptscriptstyle 0} \defi \frac{z_{\scriptscriptstyle 0}-1}{z_{\scriptscriptstyle 0}+1}=\frac{e^{j\,\theta_{\scriptscriptstyle 0}}-1}{e^{j\,\theta_{\scriptscriptstyle 0}}+1}=i\,\frac{\sin(\theta_{\scriptscriptstyle 0})}{1+\cos(\theta_{\scriptscriptstyle 0})},
\]
which shows that $z_{\scriptscriptstyle 0}$ is a pole of $G(z)$ if and only if $i\omega_{\scriptscriptstyle 0}$, with $\omega_{\scriptscriptstyle 0}\defi \frac{\sin(\theta_{\scriptscriptstyle 0})}{1+\cos(\theta_{\scriptscriptstyle 0})}>0$, is a purely imaginary pole of $G_c(s)$.
Moreover,
$i\,[G(e^{i\,\theta})-G(e^{i\,\theta})^\ast]\ge 0$ for all $\theta \in (0,\pi)$ such that $e^{i\,\theta}$ is not a pole of $G(z)$ if and only if
$i\,[G_c(i\omega)-G_c(i\omega)^\ast]\ge 0$ for all $\omega \in (0,\infty)$ such that $i\omega$ is not a pole of $G(z)$;

\item{ \hspace{-0cm} {\bf {\em (iii)}}} Let $z_{\scriptscriptstyle 0}\defi e^{i\,\theta_{\scriptscriptstyle 0}}$ with $\theta_{\scriptscriptstyle 0} \in (0,\pi)$.
Then $z_{\scriptscriptstyle 0}$, with $\theta_{\scriptscriptstyle 0} \in (0,\pi)$, is a pole of $G$ if and only if $i\omega_{\scriptscriptstyle 0}$, with $\omega_{\scriptscriptstyle 0}\defi \frac{\sin(\theta_{\scriptscriptstyle 0})}{1+\cos(\theta_{\scriptscriptstyle 0})}>0$, is a purely imaginary pole of $G_c$.
Moreover, they are poles with the same multiplicity.
 Finally,  $z_{\scriptscriptstyle 0}$ is a simple pole of $G(z)$ with residue being the matrix $K$ if and only if $i\omega_{\scriptscriptstyle 0}$ is a simple pole of $G_c(s)$ with residue being the matrix $H\defi\frac{e^{-i\,\theta_{\scriptscriptstyle 0}}}{1+\cos(\theta_{\scriptscriptstyle 0})}K$.
Notice that the normalized residual matrix $K_0$ of $G(z)$ (as defined in (\ref{normalresidu})) is positive semi-definite if and only if
$\frac{i}{z_{\scriptscriptstyle 0}}K$ is positive semi-definite and, hence, if and only if $iH$ is positive semi-definite.

%
%

 \item{ \hspace{-0cm} {\bf {\em (iv)}}} $z_{\scriptscriptstyle 0} = 1$ is a pole of $G(z)$ if and only if $s_{\scriptscriptstyle 0}=0$ is a pole of $G_c$. This fact follows straightforwardly from (\ref{use}). Moreover, they are poles with the same multiplicity. If this multiplicity is strictly greater than $2$, then $G(z)$ is trivially not D-NI.
 If this multiplicity is at most $2$, then the residual $A_{s1}$ and the quadratic residual $A_{s2}$ corresponding to $s_{\scriptscriptstyle 0}$ are related to
 the residual $A_{1}$ and the quadratic residual $A_{2}$ corresponding to $z_{\scriptscriptstyle 0}$ by: $A_{s2}=\frac{1}{4}A_2$ and $A_{s1}=\frac{1}{2}(A_1-A_2)$, since
 \[
 G(z)=G_1(z)+\frac{A_1}{z-1}+\frac{A_2}{(z-1)^2},
 \]
 where $G_1(z)$ is analytic in an open set containing $z_{\scriptscriptstyle 0}=1$, and
 \beann
 G_c(s) \ns&\ns = \ns&\ns G_{c,1}(s)+\frac{A_1}{\frac{1+s}{1-s}-1}+
 \frac{A_2}{\Big(\frac{1+s}{1-s}-1\Big)^2} \\
 \ns&\ns = \ns&\ns \Big(G_{c,1}(s)-\frac{A_1}{2}+\frac{A_2}{4}\Big)+\frac{A_1-A_2}{2\,s}+\frac{A_2}{4\,s^2} \\
 \ns&\ns = \ns&\ns \Big(G_{c,1}(s)-\frac{A_{1}}{2}+\frac{A_2}{4}\Big)+\frac{A_{s1}}{s}+\frac{A_{2s}}{s^2},
 \eeann
where $G_{c,1}(s)-\frac{A_{1}}{2}+\frac{A_2}{4}$ is analytic in an open set containing $s_{\scriptscriptstyle 0}=0$.
\item{ \hspace{-0cm} {\bf {\em (v)}}} $z_{\scriptscriptstyle 0} = -1$
is a pole of $G(z)$ if and only if $\infty$ is a pole of $G_c$. Moreover, they are poles with the same multiplicity.
If this multiplicity is strictly greater than $2$, then $G(z)$ is trivially not D-NI.
 In the case in which this multiplicity is at most $2$,
 the first coefficient $A_{s1}$ and the second coefficient $A_{s2}$ in the expansion at infinity of $G_c(s)$ are connected to
 the residual $A_{1}$ and the quadratic residual $A_{2}$ corresponding to $z_{\scriptscriptstyle 0}$ by: $A_{s2}=\frac{1}{4}A_2$ and $A_{s1}=-\frac{1}{2}(A_1+A_2)$ since
 \[
 G(z)=G_1(z)+\frac{A_1}{z+1}+\frac{A_2}{(z+1)^2},
 \]
 where $G_1(z)$ is analytic in an open set containing $z_{\scriptscriptstyle 0}=-1$, and
 \beann
 G_c(s) \ns&\ns = \ns&\ns G_{c,1}(s)+\frac{A_1}{\frac{1+s}{1-s}+1}+
 \frac{A_2}{\Big(\frac{1+s}{1-s}+1\Big)^2} \\
 \ns&\ns = \ns&\ns \Big(G_{c,1}(s)+\frac{A_1}{2}+\frac{A_2}{4}\Big)-\frac{A_1+A_2}{2}\,s+\frac{A_2}{4}\,s^2 \\
 \ns&\ns = \ns&\ns \Big(G_{c,1}(s)+\frac{A_1}{2}+\frac{A_2}{4}\Big)+A_{1s}\,s+A_{2s}\,s^2,
 \eeann
where $G_{c,1}(s)+\frac{A_1}{2}+\frac{A_2}{4}$ is rational and proper.

 \end{description}
 \end{enumerate}
Now, we apply Lemma \ref{lemmarestriction} -- see also \cite[Lemma 3.1]{Ferrante-N-13} -- in both directions and get the desired result.
\endproof

\begin{remark}\label{remsymm}
{\em
In Definition \ref{dni} we need to assume symmetry of the transfer function matrix in order to introduce the notion of a NI system as a property that is defined in the domain of analyticity: this definition is analogue to the classic definition of PR systems and has the important advantage of considering a general setting that does not require rationality assumptions.
Note, however, that if one is only interested in the rational case, it is possible to consider conditions (i)-(v) of Lemma~\ref{resd} as the definition of rational NI transfer functions and this clearly does not require any symmetry assumption. This is indeed the route taken in the first paper on continuous-time NI systems, see \cite{Lanzon-P-08,Petersen-Lanzon-10}.
A similar observation can made for the definition of strictly negative imaginary systems given below. \\
The reader can check that, as long as one considers only rational transfer functions, all the results derived in this paper can be generalized to the case of non-symmetric transfer functions.}
\end{remark}

{ We now define the notions of strictly negative imaginary systems in discrete-time.}
\begin{definition}
\label{D-SSNI}
Let $G : \complex \longrightarrow {\complex^{m \times m}}$ be a discrete-time, real transfer function. Then, $G(z)$ is {\em discrete strongly strictly negative imaginary} (D-SSNI) if
for some $\delta\in (0,1)$, the transfer function $G(\delta z)$ is D-NI and $i[G(z)-G(1/z)^\top]$ has full normal rank.
\end{definition}

\begin{remark}
{\em
The full normal rank condition is essential in the above definition as this class of systems will be needed for internal stability of positive feedback interconnections of D-NI and D-SSNI systems. If we were not to impose the full normal rank condition on the D-SSNI class, then the feedback interconnection of a D-NI system and a D-SSNI system would not be internally stable as demonstrated via the following simple example: Let $P(z)=\bsmat 1 && 1 \\[1mm] 1 && 1 \esmat$ which is clearly D-NI and let $Q(z)=\frac{2}{2\,z+1}\bsmat 1 && 1 \\[1mm] 1 && 1 \esmat$ which fulfils all properties of D-SSNI except the full normal rank condition. The positive feedback interconnection of $P(z)$ and $Q(z)$ is not internally stable as there exists a closed-loop pole at $z=3.5$.
}
\end{remark}

Now, we show that D-SSNI as defined in Definition~\ref{D-SSNI} can be equivalently checked via conditions on the domain of analyticity.
\begin{lemma}
\label{D-SSNIonOutsideDisc}
Let $G : \complex \longrightarrow {\complex^{m \times m}}$ be a discrete-time, real transfer function. Then $G(z)$ is D-SSNI if and only if there exists $\delta\in(0,1)$ such that
\begin{description}
\item{ \hspace{-0.2cm}{\bf {\em (i)}\;}} $G(z)$ is analytic in $\{z\in \complex\,:\,\,|z|>\delta\}$;
\item{ \hspace{-0.2cm}{\bf {\em (ii)}\;}} $i\,[ G(z)-G(z)^\ast ] > 0$ for all $z\in\{z\in \complex\,:\,\,|z|>\delta \; \text{and} \;\, \mathfrak{Im}(z)>0\}$;
\item{ \hspace{-0.2cm}{\bf {\em (iii)}\;}} $i\,[ G(z)-G(z)^\ast ] = 0$ for all $z\in\{z\in \complex\,:\,\,|z|>\delta \; \text{and} \;\, \mathfrak{Im}(z)=0\}$;
\item{ \hspace{-0.2cm}{\bf {\em (iv)}\;}} $i\,[ G(z)-G(z)^\ast ] < 0$ for all $z\in\{z\in \complex\,:\,\,|z|>\delta \; \text{and} \;\, \mathfrak{Im}(z)<0\}$.
\end{description}
\end{lemma}

\proof
Definition~\ref{D-SSNI} trivially gives equivalence to the existence of $\delta\in(0,1)$ such that conditions {\bf {\em (i)}}-{\bf {\em (iv)}} are satisfied with non-strict inequalities in {\bf {\em (ii)}} and {\bf {\em (iv)}} on $i\,[ G(z)-G(z)^\ast ]$. We hence only need to show that the fact that $G$ is D-SSNI implies that the inequalities in {\bf {\em (ii)}} and {\bf {\em (iv)}} are indeed strict.
The proof is similar to the one of Lemma \ref{C-SSNIonAllRHP}, the only difference being the fact that the compact set $\mathcal C$ in this case is half an annulus obtained by taking the points with non-negative imaginary parts of the annulus corresponding to the circles centered in the origin and with radii $1-\varepsilon$ and $M$, with $M$ being arbitrarily large.
\endproof

We now specialize Lemma~\ref{D-SSNIonOutsideDisc} to the unit disc. However, first we need a preliminary lemma.

\begin{lemma}
\label{C-SSNI-pre-pre-theorem-d}
Let $g : \complex \longrightarrow {\complex}$ be a scalar discrete-time, real, rational, proper transfer function. Assume that $g(z)$ is a D-SSNI function.
 If $g(1)=0$ then the multiplicity of the zero in $1$ of $g(z)$ is equal to $1$.
 Similarly, if $g(-1)=0$ then the multiplicity of the zero in $-1$ of $g(z)$ is equal to $1$.\end{lemma}
\proof
Since $g(z)$ is a D-SSNI function, it has no poles in $1$ and we can expand $g(z)$ at $1$ as
$$
g(z)=\sum_{k=h}^\infty r_k (z-1)^k,
$$
where $h$ is the multiplicity of the zero in $1$ of $g$.
Let $z=1+\varepsilon e^{i\theta}$, $0<\theta<\pi$.
If $\varepsilon$ is sufficiently small,
$
i[g(z)-g(z)^\ast] $
has the same sign of
$ -2 r_h \varepsilon^h\sin(h\theta)$
so that it can be positive for any $\theta\in(0,\pi)$ only if $h=1$.
The proof for $-1$ is similar.
\endproof

\begin{theorem}
\label{D-SSNItheorem}
Let $G : \complex \longrightarrow {\complex^{m \times m}}$ be a discrete-time, real, rational, proper transfer function. Then $G(z)$ is D-SSNI if and only if
\begin{description}
\item{ \hspace{-0.2cm}{\bf {\em (i)}\;}} $G(z)$ has all its poles with magnitude strictly less than unity;
\item{ \hspace{-0.2cm}{\bf {\em (ii)}\;}} $i\,[G(e^{i\theta})-G(e^{i\theta})^\ast]>0$ for all $\theta \in (0,\pi)$;
\item{ \hspace{-0.2cm}{\bf {\em (iii)}\;}}
$$
Q_{\scriptscriptstyle 0} \defi\lim_{\theta\rightarrow 0^+}\frac{1}{\sin(\theta)} i[G(e^{i\theta})-G(e^{i\theta})^\ast]>0
$$
\item{ \hspace{-0.2cm}{\bf {\em (iv)}\;}}
$$
Q_\pi \defi\lim_{\theta\rightarrow \pi^-}\frac{1}{\sin(\theta)} i[G(e^{i\theta})-G(e^{i\theta})^\ast]>0
$$
\end{description}
\end{theorem}

\proof
Necessity of {\bf {\em (i)}} and {\bf {\em (ii)}} is trivial from Lemma~\ref{D-SSNIonOutsideDisc}.
We now prove necessity of {\bf {\em (iii)}} (necessity of {\bf {\em (iv)}} is similar).
Assume that $G$ is D-SSNI. Then clearly the limit $Q$ defined in {\bf {\em (iii)}} exists and is positive semi-definite.
Assume by contradiction that $Q$ is singular and let $v\in\ker Q$.
Let $g'(z) \defi v\tp G(z) v$ and $g(z) \defi g'(z)-g'(1)$. Clearly, $g(z)$ is a rational proper D-SSNI function with a zero in $1$ and such that
\begin{equation}
\label{max-zero-origin-1-d}
\lim_{\theta\rightarrow 0^+}\frac{1}{\sin(\theta)} i\,[g(e^{i\theta})-g(e^{i\theta})^\ast]=0.
\end{equation}
By expanding $g(z)$ around $1$ as
$$
g(z)=\sum_{k=h}^\infty r_k (z-1)^k
$$
we see that (\ref{max-zero-origin-1-d}) implies that $h>1$, which is a contradiction in
view of Lemma \ref{C-SSNI-pre-pre-theorem-d}.

As for sufficiency,
assume that $G(s)$ is real symmetric and rational and that it satisfies
{\bf {\em (i)}}, {\bf {\em (ii)}}, {\bf {\em (iii)}} and {\bf {\em (iv)}}.
We now show that we can choose $\rho<1$ in such a way that
\begin{equation}\label{tesi-intermedia-d}
i\,[G(\rho e^{i\theta})-G(\rho e^{i\theta})^\ast]>0,\ \forall\ \theta \in (0,\pi).
\end{equation}
In view of condition {\bf {\em (ii)}}, we have that for all $\pi>\theta_2>\theta_1>0$ there exists
$\rho<1$ such that
\begin{equation}
\label{tesi-intermedia-1-d}
i\,[G(\rho e^{i\theta})-G(\rho e^{i\theta})^\ast]>0,\ \forall\ \theta \in [\theta_1,\theta_2],
\end{equation}
so that
it is sufficient to show that given an arbitrarily small $\theta_1$ and an arbitrarily large $\theta_2$, there exists
$\varepsilon>0$ such that
\begin{equation}\label{tesi-intermedia-2-d}
i\,[G(\rho e^{i\theta})-G(\rho e^{i\theta})^\ast]>0,\ \forall\ \theta \in (0,\theta_1)
\end{equation}
and
\begin{equation}\label{tesi-intermedia-3-d}
i\,[G(\rho e^{i\theta})-G(\rho e^{i\theta})^\ast]>0,\ \forall\ \theta \in (\theta_2,\pi).
\end{equation}

As for (\ref{tesi-intermedia-2-d}), let $\delta \defi\rho e^{i\theta}-1$ and consider the following expansion of $G(\delta)$:
$$
G(\delta)=D_{\scriptscriptstyle 0}+\delta D_1+\delta^2 D_2+\dots
$$
which clearly converges for $\delta$ sufficiently small (if we considered a minimal realization $G(z)=C(zI-A)^{-1}B+D$, we would have $D_{\scriptscriptstyle 0} \defi D-C(I-A)^{-1}B$ and $D_i \defi-C(I-A)^{-i-1}B$, for $i>1$).
Notice that since $G(z)$ is real symmetric by standing assumption, $D_i=D_i\tp$.
Moreover, $Q_{\scriptscriptstyle 0} \defi\lim_{\theta\rightarrow 0^+} (1/\sin(\theta))i\,[G( e^{i\theta})-G( e^{i\theta})^\ast]=-2D_1$, so that by {\bf {\em (iv)}}, we have $D_1<0$.
A direct calculation gives
$$
i\,[G(\rho e^{i\theta})-G(\rho e^{i\theta})^\ast]=-\rho\,\sin(\theta) \,2\, D_1 + i\,\sum_{j=2}^\infty [\delta^j-(\delta^\ast)^j]\,D_j.$$
Now we observe that
$$i\sum_{j=2}^\infty [\delta^j-(\delta^\ast)^j]D_j=
-2\,\rho\,\sin(\theta) \sum_{j=2}^\infty \sum_{k=0}^{j-1} [\delta^k\,(\delta^\ast)^{j-1-k}]\,D_j,$$
so that
\beann
\| i\sum_{j=3}^\infty [\delta^j-(\delta^\ast)^j]\,D_j\| \ns&\ns \leq \ns&\ns 
2\,\rho\,\sin(\theta) \sum_{j=2}^\infty j |\delta|^{j-1}\|D_j\| \\
\ns&\ns=\ns&\ns
2\,\rho\,\sin(\theta) |\delta| \sum_{j=2}^\infty j |\delta|^{j-2}\|D_j\| \\
\ns&\ns\leq\ns&\ns
2\,\rho\,\sin(\theta)|\delta| \sigma 
\eeann
for a certain $\sigma$ which remains bounded as $|\delta|$ tends to zero.
Since, by choosing a sufficiently small $\delta$ we can make
$-D_1>\sigma|\delta| I$, we have (\ref{tesi-intermedia-2}).

The proof of (\ref{tesi-intermedia-3-d}) is symmetric.
\endproof

In analogy with the continuous-time case,
we introduce the following definition of a weaker notion of strictly negative imaginary systems.

\begin{definition}
\label{defweak-d}
The discrete-time, real, rational, proper transfer function $G: \complex \longrightarrow \complex^{m \times m}$ is {\em discrete weakly strictly negative imaginary} (D-WSNI) if it satisfies conditions {\bf {\em (i)}} and {\bf {\em (ii)}} of Theorem~\ref{D-SSNItheorem}.
\end{definition}

The next lemma shows that the definition of D-WSNI characterizes properties on the outside of the unit disc too.
\begin{lemma}
\label{D-WSNIonOutsideDisc}
Let $G : \complex \longrightarrow {\complex^{m \times m}}$ be a discrete-time, real, rational, proper transfer function. Then, $G(z)$ is D-WSNI if and only if there exists $\delta\in(0,1)$ such that
\begin{description}
\item{ \hspace{-0.2cm}{\bf {\em (i)}\;}} $G(z)$ is analytic in $\{z\in \complex\,:\,\,|z|>\delta\}$;
\item{ \hspace{-0.2cm}{\bf {\em (ii)}\;}} $i\,[ G(z)-G(z)^\ast ] > 0$ for all $z\in\{z\in \complex\,:\,\,|z|\geq 1 \; \text{and} \;\, \mathfrak{Im}(z)>0\}$;
\item{ \hspace{-0.2cm}{\bf {\em (iii)}\;}} $i\,[ G(z)-G(z)^\ast ] = 0$ for all $z\in\{z\in \complex\,:\,\,|z|\geq 1 \; \text{and} \;\, \mathfrak{Im}(z)=0\}$;
\item{ \hspace{-0.2cm}{\bf {\em (iv)}\;}} $i\,[ G(z)-G(z)^\ast ] < 0$ for all $z\in\{z\in \complex\,:\,\,|z|\geq 1 \; \text{and} \;\, \mathfrak{Im}(z)<0\}$.
\end{description}
\end{lemma}
\proof
Sufficiency is trivial by restricting on $\{z\in \complex\,:\,\,|z|=1\}$. Necessity can be proven as follows: if $G$ is D-WSNI, then {\bf {\em (i)}} is satisfied and $G$ is D-NI (from Lemma \ref{resd}). If $G$ is D-NI, then {\bf {\em (ii)}}-{\bf {\em (iv)}} in Definition \ref{dni} are satisfied. Appending the $\{z\in \complex\,:\,\,|z|=1\}$ properties of $G$ to the conditions {\bf {\em (ii)}}-{\bf {\em (iv)}} in Definition \ref{dni} (since $G$ is D-WSNI) yields {\bf {\em (ii)}}-{\bf {\em (iv)}} above since $G$ fulfils {\bf {\em (i)}} above.
\endproof

{ The following lemma relates the strong classes with the weak classes with the non-strict classes of negative imaginary systems.}
\begin{lemma}
The set of D-SSNI (resp.~C-SSNI) systems is contained in the set of D-WSNI (resp.~C-WSNI) systems which is in turn contained in the set of D-NI (resp.~C-NI) systems.
\end{lemma}
\proof
Trivial from definitions.
\endproof

{ The following lemma relates a D-NI system with a D-PR system.}
\begin{lemma}
\label{lem3}
Let $G: \complex \longrightarrow \complex^{m \times m}$ be a discrete-time, real, rational, proper transfer function with no poles at $z=-1$. Then, $G(z)$ is D-NI if and only if
\bea
\label{defH}
F(z)=\frac{z-1}{z+1}\left[ G(z)-G(-1) \right]
\eea
is D-PR and $G(\infty)=G^\top(\infty)$.
\end{lemma}
\proof (Only if).
The set of poles of $F(z)$ is contained in the set of poles of $G(z)$ (in fact, in (\ref{defH}) the pole in
$-1$ of $\frac{z-1}{z+1}$ is cancelled by the zero in $-1$ of $\left[ G(z)-G(-1) \right]$). Since $G(z)$ is a symmetric, real, rational, proper, NI transfer function, $F(z)$ is analytic in $|z|>1$.   
 Let $\theta_{\scriptscriptstyle 0} \in (0,\pi)$, and assume that $z=e^{i\,\theta_{\scriptscriptstyle 0}}$ is not a pole of $G(z)$. Then, $z=e^{i\,\theta_{\scriptscriptstyle 0}}$ is not a pole of $F(z)$, and
a simple calculation gives

\beann
F(e^{i\,\theta_{\scriptscriptstyle 0}})+F(e^{i\,\theta_{\scriptscriptstyle 0}})^\ast \ns&\ns = \ns&\ns
\frac{\sin \theta_{\scriptscriptstyle 0}}{1+\cos \theta_{\scriptscriptstyle 0}}\,i\, \left[G(e^{i\,\theta_{\scriptscriptstyle 0}})-G(e^{i\,\theta_{\scriptscriptstyle 0}})^\ast\right]\ge 0
\eeann
in view of Lemma \ref{resd}.

Let us now assume that $z=e^{i\,\theta_{\scriptscriptstyle 0}}$, with $\theta_{\scriptscriptstyle 0} \in (0,\pi)$, is a pole of $G(z)$. From Lemma \ref{resd}, it is a simple pole, and from (\ref{defH}) it is also a simple pole of $H(z)$. We can write
\[
G(z)=G_1(z)+\frac{A}{z-e^{i\,\theta_{\scriptscriptstyle 0}}},
\]
where $G_1(z)$ is a rational function which is analytic in an open set containing $z=e^{i\,\theta_{\scriptscriptstyle 0}}$ and the matrix $A$ is non-zero. Then,
\beann
K_{\scriptscriptstyle 0} \ns&\ns = \ns&\ns e^{-i\,\theta_0}\,\lim_{z \to e^{i\,\theta_{\scriptscriptstyle 0}}} (z-e^{i\,\theta_{\scriptscriptstyle 0}})\,i\,G(z)\\
 \ns&\ns = \ns&\ns e^{-i\,\theta_0}\,\lim_{z \to e^{i\,\theta_{\scriptscriptstyle 0}}} (z-e^{i\,\theta_{\scriptscriptstyle 0}})\,i\,\Big(G_1(z)+\frac{A}{z-e^{i\,\theta_{\scriptscriptstyle 0}}}\Big)=i\,e^{-i\,\theta_0}\,A
\eeann
is Hermitian and positive semidefinite. The normalized residue of $F(z)$ in $e^{i\,\theta_{\scriptscriptstyle 0}}$ is given by
\beann
e^{-i\,\theta_0}\, \lim_{z\to e^{i\,\theta_{\scriptscriptstyle 0}}} (z-e^{i\,\theta_{\scriptscriptstyle 0}})\,F(z) \ns&\ns = \ns&\ns e^{-i\,\theta_0}\,
\lim_{z\to e^{i\,\theta_{\scriptscriptstyle 0}}} \frac{z-1}{z+1} \left[ (z-e^{i\,\theta_{\scriptscriptstyle 0}})\,G(z)-(z-e^{i\,\theta_{\scriptscriptstyle 0}})\,G(-1)\right] \\
\ns&\ns = \ns&\ns e^{-i\,\theta_0}\,\frac{e^{i\,\theta_{\scriptscriptstyle 0}}-1}{e^{i\,\theta_{\scriptscriptstyle 0}}+1} A=\frac{\sin \theta_{\scriptscriptstyle 0}}{1+\cos \theta_{\scriptscriptstyle 0}} i\,e^{-i\,\theta_0}\,A \ge 0.
\eeann
Let us now consider the case $\theta_{\scriptscriptstyle 0}=0$, i.e., $z_{\scriptscriptstyle 0}=e^{i\,\theta_{\scriptscriptstyle 0}}=1$. If $G(z)$ has no poles at $z=1$, neither does $F(z)$. In this case, $F(1)=0$, which gives $F(1)+F(1)^\ast=0 \ge 0$. If $G(z)$ has a simple pole at $z=1$, then $F(z)$ has no poles at $z=1$. In this case, $G(z)=G_1(z)+\frac{A}{z-1}$, where $G_1(z)$ is a rational function which is analytic in an open set containing $z=1$, and where $A \ge 0$ from {\bf {\em (iv)}} in Lemma \ref{resd} (because the quadratic residual is zero). Thus,
\beann
F(z) \ns&\ns = \ns&\ns \frac{z-1}{z+1} \left[ G_1(z)+\frac{A}{z-1}-G(-1) \right],
\eeann
so that $F(1)=A/2$, and $F(1)+F(1)^\ast=A \ge 0$. Now, consider the case in which
$G(z)$ has a double pole at $z=1$. In this case, we can write $G(z)=G_1(z)+\frac{A_1}{z-1}+\frac{A_2}{(z-1)^2}$, where $G_1(z)$ is a rational function which is analytic in an open set containing $z=1$, $A_1 \ge A_2$ and $A_2 \ge 0$. In this case,
\beann
F(z) \ns&\ns = \ns&\ns \frac{z-1}{z+1} \left[ G_1(z)+\frac{A_1}{z-1}+\frac{A_2}{(z-1)^2}-G(-1) \right] \\
\ns&\ns = \ns&\ns \left[ \frac{z-1}{z+1}\,G_1(z)+\frac{A_1}{z+1}-\frac{z-1}{z+1}\,G(-1)
-\frac{A_2}{2\,(z+1)}
\right]+\frac{A_2}{2\,(z-1)}.
\eeann
Since $G_1(z)$ is analytic in an open set containing $z=1$, such is also $\frac{z-1}{z+1}\,G_1(z)+\frac{A_1}{z+1}-\frac{z-1}{z+1}\,G(-1)
-\frac{A_2}{2\,(z+1)}$. Thus, $F(z)$ has a simple pole at $z=1$, and the corresponding residue $A_2/2$ is positive semidefinite  (in this case the residue and the normalized residue coincide because $z_{\scriptscriptstyle 0}=1$).\\
Let us finally consider the case $\theta_{\scriptscriptstyle 0}=\pi$, i.e., $z_{\scriptscriptstyle 0}=e^{i\,\theta_{\scriptscriptstyle 0}}=-1$.
We know that $G(-1)$ is finite and hence $F(-1)$ is finite as well.
Moreover, $F(e^{i\,\theta_{\scriptscriptstyle 0}})+F(e^{i\,\theta_{\scriptscriptstyle 0}})^\ast$ is positive semidefinite for all
$\theta_{\scriptscriptstyle 0} \in (0,\pi)$ that is not a pole of $G(z)$. Therefore, by continuity, we have
$F(-1)+F(-1)^\ast\geq 0$.
.

(If). Let $F$ be given by (\ref{defH}). Since $F(z)$ is symmetric, real, rational, proper, discrete-time positive real and $G(-1)=G(-1)^\top$, it is sufficient to show that
\[
G_{\scriptscriptstyle 0}(z) \defi \frac{z+1}{z-1} F(z)
\]
is D-NI because $G_{\scriptscriptstyle 0}(z)$ is D-NI if and only if $G(z)=G_{\scriptscriptstyle 0}(z)+G(-1)$ is D-NI. We observe that $G_{\scriptscriptstyle 0}(z)$ is proper, symmetric, real, rational, discrete-time and analytic in $|z|>1$. Also, $F(z)$ and $G_{\scriptscriptstyle 0}(z)$ have the same poles, with the possible exception of a pole at $z=1$. Notice that $F(z)$ does not have a pole at $z=-1$ due to its construction in (\ref{defH}).
Let $z_{\scriptscriptstyle 0}=e^{i\,\theta_{\scriptscriptstyle 0}}$ with $\theta_{\scriptscriptstyle 0}\in (0,\pi]$. Assume $z_{\scriptscriptstyle 0}$ is not a pole of $F(z)$. Then, it is not a pole of $G_{\scriptscriptstyle 0}(z)$. We find
\beann
G_{\scriptscriptstyle 0}(e^{i\,\theta_{\scriptscriptstyle 0}})\ns&\ns = \ns&\ns\frac{e^{i\,\theta_{\scriptscriptstyle 0}}+1}{e^{i\,\theta_{\scriptscriptstyle 0}}-1} \,F(e^{i\,\theta_{\scriptscriptstyle 0}}) \\
\ns&\ns = \ns&\ns -\frac{i\,\sin \theta_{\scriptscriptstyle 0}}{1-\cos \theta_{\scriptscriptstyle 0}} \,F(e^{i\,\theta_{\scriptscriptstyle 0}}),
\eeann
so that
\[
 i\,[G_{\scriptscriptstyle 0}(e^{i\,\theta_{\scriptscriptstyle 0}})-G_{\scriptscriptstyle 0}(e^{i\,\theta_{\scriptscriptstyle 0}})^\ast]=\frac{\sin \theta_{\scriptscriptstyle 0}}{1-\cos \theta_{\scriptscriptstyle 0}} \,[F(e^{i\,\theta_{\scriptscriptstyle 0}})+F(e^{i\,\theta_{\scriptscriptstyle 0}})^\ast] \geq 0,
\]
because $F(e^{i\,\theta_{\scriptscriptstyle 0}})+F(e^{i\,\theta_{\scriptscriptstyle 0}})^\ast \ge 0$. 
We now assume that $z_{\scriptscriptstyle 0}=e^{i\,\theta_{\scriptscriptstyle 0}}$ with $\theta_{\scriptscriptstyle 0}\in (0,\pi)$ is a pole of $F(z)$. Then, it is also a pole of $G_{\scriptscriptstyle 0}(z)$. Since $F(z)$ is
D-PR, $z_{\scriptscriptstyle 0}$ is a simple pole. Thus, $z_{\scriptscriptstyle 0}$ is also a simple pole of $G_{\scriptscriptstyle 0}(z)$. Moreover, the matrix $K_{\scriptscriptstyle 0}=e^{-i\,\theta_0}\,\lim_{z\to e^{i\,\theta_{\scriptscriptstyle 0}}}\,(z-e^{i\,\theta_{\scriptscriptstyle 0}}) \,F(z)$ is positive semidefinite, see Theorem \ref{restriction}.
This then implies that
\beann
e^{-i\,\theta_0}\,\lim_{z \to z_{\scriptscriptstyle 0}} (z-e^{i\,\theta_{\scriptscriptstyle 0}})\,i\,G_{\scriptscriptstyle 0}(z) \ns&\ns = \ns&\ns
e^{-i\,\theta_0}\,\lim_{z \to z_{\scriptscriptstyle 0}} i\,\frac{z+1}{z-1}\,(z-e^{i\,\theta_{\scriptscriptstyle 0}}) \,F(z) \\
\ns&\ns = \ns&\ns i\,\frac{e^{i\,\theta_{\scriptscriptstyle 0}}+1}{e^{i\,\theta_{\scriptscriptstyle 0}}-1}\,K_{\scriptscriptstyle 0} \\
\ns&\ns = \ns&\ns \frac{\sin \theta_{\scriptscriptstyle 0}}{1-\cos \theta_{\scriptscriptstyle 0}}\,K_{\scriptscriptstyle 0} \ge 0.
\eeann
When $z=1$, $F(z)$ can either have no poles or a simple pole. Assume $z=1$ is not a pole. Then, $G_{\scriptscriptstyle 0}(z)=G_1(z)+\frac{2\,F(1)}{z-1}$ where $G_1(z)$ is analytic in a region near $z=1$. Then, $K_{\scriptscriptstyle 0}=\lim_{z\to 1} (z-1)\,G_{\scriptscriptstyle 0}(z)=2\,F(1)=F(1)+F(1)^\top$ (due to $F(z)$ being symmetric),
which is non-negative in view of Theorem~\ref{restriction}. \\
Assume now that $z=1$ is a simple pole of $F(z)$. We can write $F(z)=F_1(z)+\frac{A}{z-1}$, where $F_1(z)$ is analytic near $z=1$ and $0\leq A \le 2\,F_1(1)$ (via Theorem~\ref{restriction}, since $A\geq 0$ directly from the theorem statement and $0\leq F(e^{i\theta})+F(e^{i\theta})^*=F_1(e^{i\theta})+F_1(e^{i\theta})^*-A$ implies $A \le 2\,F_1(1)$ in the limit as $\theta\rightarrow 0$ due to continuity and $F_1(1)$ being symmetric).

Hence, $G_{\scriptscriptstyle 0}(z)=\frac{z+1}{z-1}\,F(z)=\frac{z+1}{z-1}\,F_1(z)+\frac{z+1}{(z-1)^2}\,A=G_2(z)+\frac{2\,F_1(1)+A}{z-1}+\frac{2\,A}{(z-1)^2}$ where $G_2(z)$ is analytic in the neighbourhood of $z=1$. Thus, the residue and the quadratic residue are $A_1=A+2\,F_1(1)$ and $A_2=2\,A$, and the condition that ensure that $F(z)$ is D-PR now guarantees that $A_2\ge 0$ and $A_1\ge A_2$, so that $G_{\scriptscriptstyle 0}(z)$ is D-NI.
\endproof

\begin{lemma}
\label{lem2}
Let $G: \complex \longrightarrow \complex^{m \times m}$ be a discrete-time, real, rational, proper, D-NI transfer function with no poles at $z=-1$. Then
\begin{itemize}
\item $G(\infty)=G^\top(\infty)$;
\item $G(-1)$ exists and $G(-1)=G^\top(-1)$.
\end{itemize}
Furthermore, let $\bmat{c|c} A & B \\ \hline C & D \emat$ be a minimal state-space  realization of $G(z)$. Then,
\begin{itemize}
\item $C\,(I+A)^{-1}\,B=B^\top\,(I+A^\top)^{-1}\,C^\top$;
\item $F(z)=\dfrac{z-1}{z+1} \left[ G(z)-G(-1)\right]$ has a state-space  realization
\[
\bmat{c|c} A & B \\ \hline C\,(A-I)(A+I)^{-1} & C\,(A+I)^{-1}\,B \emat
\]
which is minimal when $A$ has no eigenvalues at $1$. 
\end{itemize}
\end{lemma}
\proof
Since $G(z)$ is symmetric, i.e., $G(z)=G(z)^\top$ for all $|z|>1$, we obtain $G(\infty)=G(\infty)^\top$ via a limiting argument. Since $G(z)$ has no poles at $z=-1$, it follows that $G(-1)$ exists. Now, $G(z)=G(z)^\top$ for all $|z|>1$ implies that $G(-1)=G(-1)^\top$ via continuity and a limiting argument.

From $G(-1)=G^\top(-1)$ and $D=D^\top$, it immediately follows that $C\,(I+A)^{-1}\,B=B^\top\,(I+A^\top)^{-1}\,C^\top$.

Let us now consider a state-space  realization of $H(z)$. A  realization of the transfer function matrix $\dfrac{z-1}{z+1}\,I$ is given by $\bmat{c|c} -I & I \\ \hline -2\,I & I \emat$, 
while a  realization of the term $G(z)-G(-1)=C\,(z\,I-A)^{-1}+C\,(A+I)^{-1}\,B$ is given by $\bmat{c|c} A & B \\ \hline C & C\,(A+I)^{-1}\,B \emat$. Thus, a  realization for $H(z)$ is given by

\[
\bmat{cc|c} -I & C & C\,(A+I)^{-1}\,B \\ 0 & A & B \\ \hline -2\,I & C & C\,(A+I)^{-1}\,B \emat.
\]
Changing state coordinates via
\[
T=\bmat{cc} I & C\,(I+A)^{-1} \\ 0 & I \emat
\]
yields
\beann
F(z)\ns&\ns=\ns&\ns \bmat{cc|c} -I & 0 & 0 \\ 0 & A & B \\ \hline -2\,I & C\,\left[I-2\,(I+A)^{-1}\right] & C\,(A+I)^{-1}\,B \emat \\
\ns&\ns=\ns&\ns \bmat{cc|c} -I & 0 & 0 \\ 0 & A & B \\ \hline -2\,I & C\,(A-I)(I+A)^{-1} & C\,(A+I)^{-1}\,B \emat.
\eeann
This  realization is not minimal because it is easily seen that it is not completely reachable. Eliminating the non-reachable part one obtains
\beann
F(z)\ns&\ns=\ns&\ns \bmat{c|c} A & B \\ \hline C\,(A-I)(A+I)^{-1} & C\,(A+I)^{-1}\,B \emat,
\eeann
which is minimal if $\det (A-I) \neq 0$.
\endproof

\begin{remark}
{\em
Notice that we have derived condition $G(-1)= G(-1)^\top$ as a consequence of the symmetry 
of $G(z)$. However, if we consider, in the spirit of Remark \ref{remsymm}, the possibly non-symmetric case, then condition $G(-1)= G(-1)^\top$ still holds. More precisely,  assuming that rational NI systems are defined by  conditions (i)-(v) of Lemma~\ref{resd} (and that symmetry is not assumed), we have that if $-1$ is not a pole of $G(z)$ then $G(-1)= G(-1)^\top$.
In fact, since by condition (ii) of Lemma~\ref{resd},
$i\,[G(e^{i\theta})-G(e^{i\theta})^\ast]\ge 0$ for all $\theta \in (0,\pi)$ (except for the values of $\theta$ for which $z=e^{i\,\theta}$ is a pole of $G(z)$), we can use continuity and conclude that $i\,[G(-1)-G(-1)^\ast]\ge 0$, but $G(-1)$ is real so that we get that
$i\,[G(-1)-G(-1)^\top]$ is positive semi-definite. But the diagonal entries of $i\,[G(-1)-G(-1)^\top]$ are zero so that we necessarily have
$G(-1)-G(-1)^\top=0$.
Similarly, assuming that rational NI systems are defined by  conditions (i)-(v) of Lemma~\ref{resd} (and that symmetry is not assumed), we have that if $1$ is not a pole of $G(z)$ then $G(1)= G(1)^\top$.
\\
In this non-symmetric setting, it is easy to check that the result analogue to Lemma \ref{lem3} is that $G(z)$ without poles in $-1$ is NI if and only if $F(z)$ defined by (\ref{defH}) is PR {\em and 
$G(-1)= G(-1)^\top$.}
}
\end{remark}

{ We are now in a position to give a discrete-time negative imaginary lemma that gives a complete state-space characterization of D-NI systems.}

\begin{theorem}
\label{the2}
Let $\bmat{c|c} A & B \\ \hline C & D \emat$ be a minimal state-space  realization of a discrete-time, real, rational, proper transfer function $G(z)$. Suppose $\det (I+A) \neq 0$ and $\det (I-A) \neq 0$. Then, $G(z)$ is D-NI if and only if
$D=D^\top$ and there exists a real matrix $X=X^\top>0$ such that
\begin{equation}
\label{D-NI_lemma}
X-A^\top X\,A \ge 0 \quad \text{and} \quad C=-B^\top (A^\top-I)^{-1} X\,(A+I).
\end{equation}
\end{theorem}
\proof
First, note that
\bea
\label{eq1}
A\,(A-I)^{-1}=I+(A-I)^{-1}.
\eea
Now, in view of Lemma \ref{lem3}, $G(z)$ is D-NI if and only if $H(z)=\frac{z-1}{z+1}\left[ G(z)-G(-1) \right]$ is D-PR and $D=D^\top$.
By Lemma \ref{lem2}, this is equivalent to 
\[
\bmat{c|c} A & B \\ \hline C\,(A-I)(A+I)^{-1} & C\,(A+I)^{-1}\,B \emat
\]
 being D-PR and $D=D^\top$. Using Lemma \ref{theoremAnderson}, the latter conditions are equivalent to $D=D^\top$ and there exists $X>0$ and $L,W$ such that
\bea
&&X-A^\top\,X\,A=L^\top\,L \label{eqA1} \\
&&(A^\top+I)^{-1}(A^\top-I)\,C^\top -A^\top\,X\,B=L^\top\,W \label{eqA2} \\
&&C\,(A+I)^{-1}\,B+B^\top\,(I+A^\top)^{-1}\,C^\top-B^\top\,X\,B=W^\top\,W \label{eqA3}
\eea
Eq. (\ref{eqA2}) can be written as
\[
C=(W^\top L+B^\top X\,A)(A-I)^{-1} (A+I),
\]
which can be substituted into (\ref{eqA3}) to give
\beann
&&B^\top X\,[ I+(A-I)^{-1}] \,B +
B^\top \,[ I+(A^\top-I)^{-1}] \,X\,B-B^\top X\,B \\
&&\qquad \quad =W^\top W -W^\top L\,(A-I)^{-1} B-B^\top (A^\top-I)^{-1} L^\top W
\eeann
in view of (\ref{eq1}). This equation can also be written as
\beann
&&B^\top X \,(A-I)^{-1} B+B^\top (A^\top-I)^{-1} X \,B+B^\top X\,B+
B^\top (A^\top-I)^{-1}L^\top L (A-I)^{-1} B \\
&&\qquad \quad = [W-L(A-I)^{-1} B]^\top[W-L(A-I)^{-1} B].
\eeann
Plugging the term $L^\top L$ of (\ref{eqA1}) into the latter yields
\beann
&&B^\top X \,(A-I)^{-1} B+B^\top (A^\top-I)^{-1} X \,B+B^\top X\,B+
B^\top (A^\top-I)^{-1}X \,(A-I)^{-1} B\\
&& - B^\top (A^\top-I)^{-1}A^\top\,X\,A \,(A-I)^{-1} B \\
&&\qquad \quad = [W-L(A-I)^{-1} B]^\top[W-L(A-I)^{-1} B].
\eeann
Using (\ref{eq1}), it is easily seen that the left hand-side of this equation is equal to zero, so that $W=L\,(A-I)^{-1} B$. This means that $G(z)$ is D-NI if and only if $D=D^\top$,
\begin{itemize}
\item $X-A^\top X\,A \ge 0$;
\item $C=\left(B^\top (A^\top-I)^{-1} (X-A^\top X \,A)+B^\top X\,A\right) (A-I)^{-1}(A+I)$.
\end{itemize}
Now, using (\ref{eq1}), $G(z)$ is D-NI if and only if $D=D^\top$ and there exists $X>0$ such that
\begin{itemize}
\item $X-A^\top X\,A \ge 0$
\item $C\,(A+I)^{-1} =-B^\top (A^\top-I)^{-1} X$
\end{itemize}
\endproof

\begin{remark}
{\em
In Theorem~\ref{the2}, note that given matrices $A, B, C$,
\begin{align*}
&\exists X=X^\top>0: X-A^\top X\,A \ge 0 \quad \text{and} \quad C=-B^\top (A^\top-I)^{-1} X\,(A+I) \\
\Leftrightarrow\quad &\exists Y=Y^\top>0: Y-A Y\,A^\top \ge 0 \quad \text{and} \quad B=-(A-I) Y\,(A^\top+I)^{-1} C^\top.
\end{align*}
}
\end{remark}

Similar to Lemma 2 in \cite{Lanzon-P-08}, we here show that the gain of the system at $G(1)$ and $G(-1)$ can be ordered as given in the following lemma.

\begin{lemma}
Let $G: \complex \longrightarrow \complex^{m \times m}$ be a discrete-time, real, rational, proper, D-NI (resp.~D-WSNI) transfer function with no poles at $+1$ and $-1$. Then
\[
G(1)-G(-1) \ge 0 \;(\text{resp.} > 0).
\]
\end{lemma}
\proof Using Theorem \ref{the2} and a minimal  realization for $G(z)$, we find
\beann
G(1)-G(-1) \ns&\ns = \ns&\ns C\,(I-A)^{-1}\,B+D-C(-I-A)^{-1}\,B-D \\
\ns&\ns = \ns&\ns C\,\left[ (I-A)^{-1}+(I+A)^{-1} \right]\,B \\
\ns&\ns = \ns&\ns 2\,C\,(I+A)^{-1} (I-A)^{-1}\,B \\
\ns&\ns = \ns&\ns -2\,B^\top (A^\top-I)^{-1}\,X\,(I-A)^{-1}\,B \\
\ns&\ns = \ns&\ns 2\,B^\top (I-A)^{-\top} X (I-A)^{-1}\,B \ge 0.
\eeann
This concludes the proof for $G$ being D-NI. \\

Now, we focus on $G$ being D-WSNI. The strict inequality result will be proven via a contra-positive argument. Suppose there exists an $x\in\real^{m}$ such that $[G(1)-G(-1)]x=0$. Then $B^\top (I-A)^{-\top} X (I-A)^{-1}\,Bx=0$ which implies that $Bx=0$ as $X>0$. This then implies that $G(e^{i\theta})x=Dx\;\;\forall\theta\in(0,\pi)$, {i.e.
$$
(G(e^{i\theta})-D)x=0,\quad \forall\theta\in(0,\pi).
$$
But since $G$ is WSNI, $i\,[G(e^{i\theta})-G(e^{i\theta})^*]$ is positive definite for all $\theta\in(0,\pi)$ so that if, for
$\theta_0\in(0,\pi)$, $x$ is such that  $x^*[i\,[G(e^{i\theta_0})-G(e^{i\theta_0})^*]]x=0$, we can conclude that $x=0$.
Now recall that $D=D^\top$. Hence,
$x^*[i\,[G(e^{i\theta_0})-G(e^{i\theta_0})^*]]x=i\, x^*[(G(e^{i\theta_0})-D)-(G(e^{i\theta_0})-D)^*]x=0$.
Hence $x=0$, so that $[G(1)-G(-1)]$} must be nonsingular. This completes the proof.
\endproof

{ The following result shows under what circumstances are D-NI, D-WSNI and D-SSNI properties preserved when such systems are interconnected in feedback.}
\begin{lemma}\label{Redheffer}
Let $S_1:\complex\rightarrow\complex^{m_1\times m_1}$ be D-NI (resp. D-WSNI or D-SSNI) and $S_2:\complex\rightarrow\complex^{m_2\times m_2}$ be D-NI (resp.~D-WSNI or D-SSNI). Let $0<a,b\leq\min\{m_1,m_2\}$
and suppose the feedback interconnection corresponding to the Redheffer Star product $S_1\star S_2$ be internally stable.\footnote{This is the standard meaning of ``internal stability'',
i.e.~add two extra exogenous input signals to the internal signals and ensure that all output signals and all internal signals are energy-bounded for any energy-bounded exogenous input excitation.}
Then $S_1\star S_2$ is D-NI (resp. D-WSNI or D-SSNI).

Furthermore, if
\begin{itemize}
\item $a=b=m_2$, then $S_1\star S_2 = F_l(S_1,S_2)$;
\item $a=b=m_1$, then $S_1\star S_2 = F_u(S_2,S_1)$;
\item $a=b=m_2$, $S_1=\begin{bmatrix} P & I_a \\ I_a & 0\end{bmatrix}$ and $S_2=Q$, then $S_1\star S_2 = P+Q$;
\item $a=b=m_1/2=m_2/2$, $S_1=\begin{bmatrix} 0 & I_a\\ I_a & P\end{bmatrix}$ and $S_2=\begin{bmatrix} Q & I_a \\ I_a & 0\end{bmatrix}$, then
$S_1\star S_2 = \begin{bmatrix}-P & I_a \\ I_a & -Q\end{bmatrix}^{-1}$ $=\begin{bmatrix}Q(I_a-PQ)^{-1} & (I_a-QP)^{-1} \\ (I_a-PQ)^{-1} & P(I_a-QP)^{-1}\end{bmatrix}$
which corresponds to the positive feedback interconnection $[P,Q]$.
\end{itemize}
\end{lemma}
\proof
Given $S_1(z), S_2(z)$ and complex vectors $y_1, y_2, u_1, u_2, \alpha, \beta$ of compatible dimension satisfying
$\begin{bmatrix}y_1 \\ \alpha\end{bmatrix}=S_1(z) \begin{bmatrix}u_1\\ \beta\end{bmatrix}$ and $\begin{bmatrix}\beta \\ y_2\end{bmatrix}=S_2(z) \begin{bmatrix}\alpha \\ u_2\end{bmatrix}$,
it follows that $\begin{bmatrix}y_1\\ y_2\end{bmatrix}=S_1(z)\star S_2(z)\begin{bmatrix}u_1 \\ u_2\end{bmatrix}$.
Then, for all $\begin{bmatrix}u_1\\\beta\end{bmatrix}\in\complex^{m_1}, \begin{bmatrix}\alpha\\u_2\end{bmatrix}\in\complex^{m_2}$:
\begin{align*}
 \begin{bmatrix}u_1^*&u_2^*\end{bmatrix}\left[i([S_1(z)\star S_2(z)] - [S_1(z)\star S_2(z)]^*)\right]\begin{bmatrix}u_1\\u_2\end{bmatrix} \hspace{-15em}& \\
 &= i\begin{bmatrix}u_1^*&u_2^*\end{bmatrix}\begin{bmatrix}y_1\\y_2\end{bmatrix} - i\begin{bmatrix}y_1^*&y_2^*\end{bmatrix}\begin{bmatrix}u_1\\u_2\end{bmatrix} \\
 &= i\begin{bmatrix}u_1^*&\beta^*\end{bmatrix}\begin{bmatrix}y_1\\\alpha\end{bmatrix} - i\begin{bmatrix}y_1^*&\alpha^*\end{bmatrix}\begin{bmatrix}u_1\\\beta\end{bmatrix}
 + i\begin{bmatrix}\alpha^*&u_2^*\end{bmatrix}\begin{bmatrix}\beta\\y_2\end{bmatrix} - i\begin{bmatrix}\beta^*&y_2^*\end{bmatrix}\begin{bmatrix}\alpha\\u_2\end{bmatrix} \\
 &= \begin{bmatrix}u_1^*&\beta^*\end{bmatrix}\left[i(S_1(z) - S_1(z)^*)\right]\begin{bmatrix}u_1\\\beta\end{bmatrix}
 + \begin{bmatrix}\alpha^*&u_2^*\end{bmatrix}\left[i(S_2(z) - S_2(z)^*)\right]\begin{bmatrix}\alpha\\u_2\end{bmatrix}.
\end{align*}
Since the Redheffer star interconnection is internally stable, the three respective results (D-NI, D-WSNI, D-SSNI) then follow by applying Definition~\ref{dni}, Lemma~\ref{D-SSNIonOutsideDisc}
or Lemma~\ref{D-WSNIonOutsideDisc} respectively on the corresponding domains of $z\in\complex$ for $S_1(z)$ and $S_2(z)$.

The four cases where $a, b, S_1$ and $S_2$ are restricted are trivial consequences of a Redheffer calculation.
\endproof

\begin{remark}
{\em
Lemma~\ref{Redheffer} holds also in continuous-time with all of D-NI, D-SSNI and D-WSNI replaced by C-NI, C-SSNI and C-WSNI respectively.\footnote{see \cite{Petersen-Lanzon-10} for a sub-class.}}
\end{remark}
\begin{example}
{\em
This example shows that it is not possible to mix and match properties of $S_1$ and $S_2$ for the strict results in Lemma~\ref{Redheffer} to hold.

Let $S_1=\begin{bmatrix}\begin{pmatrix}1&0\\0&1\end{pmatrix} & \begin{pmatrix}1\\0\end{pmatrix} \\ \begin{pmatrix}1&0\end{pmatrix} & 0\end{bmatrix}$ which is clearly D-NI and let $S_2=z^{-1}$ which is clearly D-SSNI (and hence also D-WSNI and hence also D-NI). Then $S_1 \star S_2=\begin{bmatrix}1+z^{-1} & 0 \\ 0 & 1\end{bmatrix}$ which is only D-NI (and not D-WSNI nor D-SSNI).
}
\end{example}

The following stability theorem here applies only to real, rational, proper systems but invokes only the interconnection of D-NI and D-WSNI systems.
{ It is the discrete-time analogue of Theorem~5 in \cite{Lanzon-P-08}.}

\begin{theorem}
 Let $P:\complex\rightarrow\complex^{m\times m}$ be a real, rational, proper, D-NI system with no poles at $+1$ and $-1$, and let $Q:\complex\rightarrow\complex^{m\times m}$ be a real, rational, proper, D-WSNI system. Suppose $P(-1)Q(-1)= 0$ and $Q(-1)\geq 0$. Then
 \[ [P,Q] \text{ is internally stable}\qquad\Leftrightarrow\qquad \bar{\lambda}\left(P(1)Q(1)\right) < 1. \]
\end{theorem}
\proof
The proof trivially follows by applying \cite[Theorem~5]{Lanzon-P-08} or \cite[Theorem~1]{Xiong-PL-10} on the systems $M(s)=P(\frac{1+s}{1-s})$ and $N(s)=Q(\frac{1+s}{1-s})$ obtained through the bilinear transformation $z=\frac{1+s}{1-s}$.
\endproof

\section*{Concluding remarks}

In this paper we presented a definition of negative imaginary systems for discrete-time systems that hinges entirely on properties of the transfer function matrix and not on a real, rational, proper, finite-dimensional realization. 
We have drawn a full picture which illustrates the relationship that exists between the notions of positive real and negative imaginary systems, as well as strictly positive real and strictly negative imaginary systems, both in continuous time and in discrete time. 
{
Indeed, notice that, pretty much as it happened for the classical theory of positive real systems, even for negative imaginary systems our definitions for the discrete-time and continuous-time cases can be viewed as a single definition referred to different analyticity domains.
In fact, we can define a function $G : \complex \longrightarrow {\complex^{m \times m}}$ analytic in an open subset $\Omega\subset\complex$, to be is {\em skew-imaginary} if{ 
\begin{itemize}
\item $i\,[ G(s)-G(s)^\ast ] \ge 0$ for all $s \in \Omega$ such that $\mathfrak{Im}\{s\} >0$;
\item $i\,[ G(s)-G(s)^\ast ] = 0$ for all $s \in \Omega$ such that $\mathfrak{Im}\{s\}=0$;
\item $i\,[ G(s)-G(s)^\ast ] \le 0$ for all $s \in \Omega$ such that $\mathfrak{Im}\{s\} <0$.
\end{itemize}}
Then, it is clear that a function is NI if it is analytic in $\Omega$ and skew-imaginary there.
Here, $\Omega$ is the open right half complex plane for the continuous case and the set $\{z\in\complex:\ |z|>1\}$ for the discrete-time case.

Finally, we have derived a stability analysis result for the interconnections of D-NI and D-WSNI systems.}

\end{document}